\documentclass[a4paper,reqno]{amsart}

\usepackage{amsmath,amssymb,amsthm}
\usepackage{fix-cm}
\usepackage{booktabs}
\usepackage{graphicx}
\graphicspath{{explicit-finite-horizon/figures/}{../figures/}}
\usepackage[
    colorlinks=true,
    linkcolor=red,
    citecolor=red,
    urlcolor=red,
    pagebackref=true
]{hyperref}

\newtheorem{theorem}{Theorem}[section]
\newtheorem{corollary}[theorem]{Corollary}
\newtheorem{lemma}[theorem]{Lemma}
\newtheorem{proposition}[theorem]{Proposition}
\newtheorem{remark}[theorem]{Remark}
\newtheorem*{claim}{Claim}

\newcommand{\R}{\mathbb R}

\newcommand{\Z}{\mathbb Z}

\newcommand{\vect}[1]{\boldsymbol{#1}}       
\newcommand{\planarvec}[1]{\boldsymbol{#1}} 
\newcommand{\matr}[1]{\mathbf{#1}}

\newcommand{\latt}[1]{\mathcal{#1}}
\newcommand{\funddom}[1]{\mathcal{F}_{#1}}
\newcommand{\norm}[1]{\left\lVert #1\right\rVert}
\newcommand{\onenorm}[1]{\left\lVert #1\right\rVert_1}
\newcommand{\abs}[1]{\left|#1\right|}
\newcommand{\dotprod}{\mathbin{\cdot}}
\DeclareMathOperator{\covol}{covol}
\DeclareMathOperator{\dist}{dist}

\DeclareMathOperator{\Span}{span}

\newcommand{\FCC}{\Lambda}
\newcommand{\BlockingMargin}{\eta_0}
\newcommand{\dualFCC}{\FCC^*}
\newcommand{\projlat}[1]{\Gamma_{\vect{#1}}}
\newcommand{\orthcomp}[1]{\vect{#1}^{\perp}}
\newcommand{\orthproj}[1]{\pi_{\orthcomp{#1}}}
\newcommand{\CL}{\mathcal C_L}
\newcommand{\CO}{\mathcal C_O}
\newcommand{\CT}{\mathcal C_T}
\newcommand{\CB}{\mathcal C_B}
\newcommand{\Horizon}{\mathcal H}
\newcommand{\Area}{\mathcal A}
\newcommand{\PeriodTorus}{\mathbb T}

\title[An Explicit Three-Dimensional Finite-Horizon Billiard]
{An Explicit Three-Dimensional Dispersing Billiard with Finite Horizon}
\author{Peyman Eslami}
\date{November 13, 2026}
\thanks{The author used OpenAI Codex to assist with editing the manuscript,
    preparing figures, and running computational checks.}

\begin{document}
    \begin{abstract}
        We construct an explicit three-dimensional periodic Lorentz gas whose
        quotient by the period lattice \(2\Z^3\) has 38 spherical scatterers
        with rational centers and four rational radii. The scatterers are
        strictly separated, so the associated dispersing billiard has no corners
        or cusps. We prove finite horizon by establishing a common blocking
        margin for every rational direction. A projected-lattice estimate for
        the face-centered cubic core leaves nine exceptional directions,
        handled by a finite line-space covering certificate. Simultaneous
        rational approximation transfers this margin to every direction and
        yields an explicit upper horizon bound. All finite calculations use exact
        arithmetic, with verification data and algorithms included. The positive
        separation and blocking margins also ensure stability under sufficiently
        small independent perturbations of the centers and radii.
    \end{abstract}

    \maketitle

    \section{Introduction}

    In a periodic Lorentz gas, a point particle moves freely among a periodic
    family of convex scatterers and reflects elastically at their boundaries.
    Taking the quotient by the period lattice gives a dispersing billiard on
    a flat torus. Finite horizon means that the lengths of all collision-free
    segments have a common upper bound. We construct a three-dimensional
    billiard with finite horizon from strictly separated spherical scatterers
    with rational centers and radii.

    \begin{theorem}[Explicit finite-horizon configuration]
        \label{thm:explicit}
        The family \(\mathcal B\) specified in Section~\ref{sec:explicit} has full
        translation lattice \(2\Z^3\) and 38 scatterers modulo this lattice. The
        scatterers are strictly separated balls with rational centers and four
        rational radii. Their complement is a three-dimensional dispersing
        billiard with no corners or cusps, and its horizon satisfies
        \[
            \Horizon<453601.
        \]
    \end{theorem}

    The strict separation and blocking inequalities have positive margins.
    Consequently, uniformly small independent perturbations of all balls in the
    periodic lift preserve strict separation and finite horizon, even when the
    perturbations destroy periodicity, as shown in
    Corollary~\ref{cor:stability}. The result concerns the geometric
    finite-horizon problem. We do not address
    further dynamical properties of this billiard.

    Existence and explicitness are distinct in this problem. Classical results
    summarized by Henk and Zong show that every lattice
    packing of balls in dimension at least three has an infinite cylindrical
    hole~\cite[p.~31]{HZ}. Thus one ball orbit under a period lattice cannot
    produce finite horizon. Henk and Zong then proved that the obstruction
    disappears for a general periodic set of centers: for every dimension
    \(n\) and every \(\varepsilon>0\), they construct a periodic point set
    \(P(n,\varepsilon)\) such that the unit balls centered at
    \(P(n,\varepsilon)\) form a packing and all segments avoiding the interiors
    of the \(\varepsilon\)-balls have uniformly bounded length~\cite[Theorem~1]{HZ}.
    Taking \(0<\varepsilon<1\) and using the \(\varepsilon\)-balls as
    scatterers gives congruent, strictly separated finite-horizon configurations.

    Sanders emphasized the remaining problem of producing an explicit model in
    three dimensions. He analyzed a simple cubic model with two radii in which
    finite horizon is obtained by allowing the spheres to
    overlap~\cite{Sanders08}. Dettmann subsequently isolated a stronger
    question: find an explicit periodic arrangement of equal-sized,
    nonoverlapping spheres with finite horizon in dimension at least
    three~\cite[Open Problem~2]{Dettmann14}.

    To our knowledge, Theorem~\ref{thm:explicit} gives the first periodic
    finite-horizon configuration of spherical scatterers in three dimensions
    for which the period lattice, finitely many center representatives, and
    radii are explicitly specified and the scatterers are strictly separated.
    The construction uses four radii and therefore does not solve Dettmann's
    equal-radius problem.
    Strict separation is stronger than nonoverlap: it rules out both corner
    intersections caused by overlap and cusp contacts caused by tangency.

    Throughout, a direction is \emph{rational} if it is parallel to a nonzero
    vector of \(\Z^3\).

    The construction has a lattice core and a finite correction. Large balls
    centered at the face-centered cubic (FCC) lattice block all but finitely
    many rational directions. Smaller balls placed in the remaining gaps block
    the exceptional directions.

    The argument has four steps:
    \begin{enumerate}
        \item Specify the rational centers and radii and verify strict
        separation (Section~\ref{sec:explicit}).
        \item Use the FCC line estimate (Lemma~\ref{lem:fcc-lines}) to reduce the
        proof of blocking in every rational direction to checking 61 unoriented
        directions. Exact planar covering radii establish blocking in 52 of
        these directions
        (Section~\ref{sec:rational-blocking}).
        \item Certify the remaining nine directions by a finite covering of
        line space. Both the core estimates and this certificate use balls
        whose radii are reduced by \(1/250\), giving a common blocking margin.
        \item Approximate arbitrary directions by rational ones. Periodicity bounds
        the distance to an intersection along the approximating line, and the
        common blocking margin allows this intersection to persist when the
        direction changes. Controlling the period and the approximation error
        yields the upper horizon bound
        (Section~\ref{sec:quantitative}).
    \end{enumerate}
    The positive blocking margin is the link between the finite covering
    calculations and the finite-horizon conclusion.

    Section~\ref{sec:explicit} begins by fixing the billiard terminology. On a
    first reading, the main argument can be followed in order, using the
    lattice facts and covering formulas proved in
    Appendix~\ref{app:lattice-facts}. The finite verification is also
    self-contained: Appendix~\ref{app:certificate} derives the
    cell-containment inequality, and Appendix~\ref{app:reproducibility}
    gives deterministic exact-arithmetic verification algorithms.

    \section{The explicit billiard table}
    \label{sec:explicit}

    Let \(\latt{L}\subset\R^3\) be a full-rank lattice. Translating a
    finite family of representative balls, with their assigned radii, by
    every vector in \(\latt{L}\) gives its \emph{periodic lift}. The
    translates of one representative form its \emph{ball orbit} under
    \(\latt{L}\). We call \(\latt{L}\) a \emph{period lattice} of the
    resulting configuration. Its \emph{full translation lattice} consists
    of all vectors whose translations preserve the entire family of balls,
    including their radii; it may be larger than the specified period lattice.

    A finite family of balls
    modulo \(\latt{L}\), with pairwise disjoint interiors in the periodic lift,
    defines a dispersing billiard on the flat torus
    \(\R^3/\latt{L}\)~\cite{CM06,Dettmann14}. The construction below uses
    \(\latt{L}=2\Z^3\), with representatives in the fundamental domain
    \(\funddom{2\Z^3}=[0,2)^3\). From now on, the ambient torus is
    \(\PeriodTorus=\R^3/2\Z^3\).

    A \emph{line} means the full image of an affine line, in \(\R^3\) or in the
    quotient torus. A line is \emph{blocked} only when it enters the interior of
    a scatterer; tangential contact does not count.
    A family of balls \emph{blocks a direction}
    if every affine line parallel to that direction enters the interior of
    at least one ball in the family. A segment is
    \emph{collision-free} if it avoids the interiors of all scatterers;
    tangential boundary contact is permitted. The
    \emph{horizon} \(\Horizon\) is the supremum of the lengths of all
    collision-free segments, and the billiard has \emph{finite horizon} when
    \(\Horizon<\infty\). A free path between successive elastic reflections is
    one such segment. A family of balls is a
    \emph{packing}, or is \emph{nonoverlapping}, if its members have pairwise
    disjoint interiors; tangency is permitted. A packing is
    \emph{strictly separated} if distinct balls
    have positive distance. In the configuration below, strict separation makes
    every boundary component a disjoint smooth sphere and therefore excludes
    corners and cusps.

    \subsection{Centers, radii, and the core construction}

    The large balls form a face-centered cubic (FCC) core, whose balls
    block all but nine rational directions even after their radii are
    reduced by \(1/250\). The three smaller-ball
    families supplement this core: together, the four families give the
    line-space coverings certified in Lemma~\ref{lem:certificate} for the
    remaining directions. We first specify their centers and radii, then
    verify separation and blocking.

    All centers below are listed modulo \(2\Z^3\). The first three
    families consist of large balls and balls at
    selected octahedral and tetrahedral interstitial
    sites of the FCC lattice: the nearest FCC centers to such a site form
    a regular octahedron or tetrahedron, respectively. Denote the corresponding center sets by
    \(\CL\), \(\CO\), and \(\CT\):

    \begin{align*}
        \CL&=\{\vect{u}\in\{0,1\}^3:u_1+u_2+u_3\text{ is even}\},\\
        \CO&=\{(1,0,0),(0,1,0),(0,0,1)\},\\
        \CT&=\left\{(1,1,1)+\frac12\vect{\varepsilon}:
                 \varepsilon_i\in\{-1,1\},\
                 \varepsilon_1\varepsilon_2\varepsilon_3=1\right\}.
    \end{align*}

    The four centers in \(\CL\) carry the large balls of radius \(7/10\). The
    three centers in \(\CO\) carry balls of radius \(7/25\), and the four
    centers in \(\CT\) carry balls of radius \(7/50\).

    The periodic lift of \(\CL\) is the FCC lattice
    \begin{equation}
        \CL+2\Z^3
        =\FCC
        =\{\vect{z}\in\Z^3:z_1+z_2+z_3\text{ is even}\}.
        \label{eq:fcc}
    \end{equation}
    The parity condition in~\eqref{eq:fcc} is the standard
    description of the FCC lattice~\cite[Chap.~4, \S6.3, p.~112]{CS99}.
    The radius-\(7/10\) balls centered at \(\FCC\) form the FCC core.
    Section~\ref{sec:rational-blocking} proves that the FCC core blocks all
    but finitely many rational directions.

    To complete the configuration, add 27 balls of radius
    \(2/25\), with center set \(\CB\). We list \(10\CB\) to give the
    coordinates as integers:
    \[
        10\CB=\left\{
                  \begin{array}{ccc}
                      (2,17,13)&(18,13,17)&(18,13,3)\\
                  (8,7,13)&(12,17,17)&(12,3,17)\\
                  (13,2,3)&(7,2,3)&(17,18,7)\\
                  (3,18,7)&(3,2,13)&(13,13,12)\\
                  (7,13,12)&(13,8,13)&(7,3,18)\\
                  (17,7,2)&(3,7,2)&(3,13,18)\\
                  (13,12,7)&(7,12,7)&(17,7,18)\\
                  (17,18,13)&(8,17,17)&(18,3,13)\\
                  (8,17,3)&(2,13,3)&(18,3,7)
                  \end{array}
        \right\}.
    \]
    Let \(\mathcal B\) denote the periodic family of closed balls obtained
    by translating the listed balls, with their assigned radii, by all
    vectors in \(2\Z^3\). There are therefore \(4+3+4+27=38\) listed scatterers on
    \(\PeriodTorus\).

    We call the three smaller-ball families represented by
    \(\CO\), \(\CT\), and \(\CB\), together with all their
    \(2\Z^3\)-translates, the \emph{correction family}. The line-space
    certificate uses this family jointly with the FCC core.

    \subsection{Separation and the period lattice}

    We now verify that the balls in \(\mathcal B\) are strictly separated
    and that \(2\Z^3\) is its full translation lattice.

    \begin{proposition}[Separation and periods]
        \label{prop:packing-period}
        The family \(\mathcal B\) is strictly separated and has full
        translation lattice \(2\Z^3\).
    \end{proposition}

    \begin{proof}
        For separation, translate the pair so that one center is a listed
        representative \(\vect{p}\in[0,2)^3\). The other has the form
        \(\vect{q}+2\vect{k}\), where \(\vect{q}\in[0,2)^3\) is listed and
        \(\vect{k}\in\Z^3\). If \(\abs{k_j}\ge2\) for some \(j\), then
        \[
            \abs{p_j-q_j-2k_j}>2.
        \]
        Since the sum of two radii is at most \(7/5<2\), such a pair is
        separated. Thus it is enough to consider
        \(\vect{k}\in\{-1,0,1\}^3\), omitting the coincident choice
        \(\vect{p}=\vect{q}\), \(\vect{k}=\vect{0}\).

        For the remaining finite set, consider the squared separation surplus
        \[
            \norm{\vect{p}-(\vect{q}+2\vect{k})}^2
            -(r_{\vect{p}}+r_{\vect{q}})^2.
        \]
        This surplus is a difference of squared
        quantities, not the distance between the two ball boundaries.
        All these numbers are rational, and the exact packing algorithm in
        Appendix~\ref{alg:packing} gives
        \[
            \min_{\vect{p},\vect{q},\vect{k}}
            \left\{\norm{\vect{p}-(\vect{q}+2\vect{k})}^2
            -(r_{\vect{p}}+r_{\vect{q}})^2\right\}
            =\frac{29}{2500}>0.
        \]
        Hence the balls are strictly separated. The minimum occurs between the
        radius-\(7/10\)
        ball at the origin and the translate by \((0,-2,-2)\) of the
        radius-\(2/25\) ball centered at \((1/5,17/10,13/10)\).
        The displacement from the origin to this translated center is
        \((1/5,-3/10,-7/10)\), and
        \[
            \norm{(1/5,-3/10,-7/10)}^2=\frac{31}{50},\qquad
            \left(\frac7{10}+\frac2{25}\right)^2=\frac{1521}{2500},
        \]
        whose difference is \(29/2500\).

        To identify the full translation lattice, suppose translation by
        \(\vect{t}\) preserves \(\mathcal B\). It sends the radius-\(7/10\)
        ball at the origin to a ball of the same radius centered at
        \(\vect{t}\). The centers of these balls form \(\FCC\)
        by~\eqref{eq:fcc}, so \(\vect{t}\in\FCC\).
        Modulo \(2\Z^3\), the four classes in
        \(\FCC\) are represented by
        \[
            (0,0,0),\qquad (0,1,1),\qquad (1,0,1),\qquad (1,1,0).
        \]
        The three nonzero classes do not preserve the radius-\(7/25\) centers. Indeed,
        \begin{align*}
            (1,0,0)+(0,1,1)&=(1,1,1),\\
        (0,1,0)+(1,0,1)&=(1,1,1),\\
        (0,0,1)+(1,1,0)&=(1,1,1),
        \end{align*}
        and \((1,1,1)\), modulo \(2\Z^3\), is not a radius-\(7/25\) center. Hence
        no nonzero class is a period, and the full translation lattice is exactly
        \(2\Z^3\).

    \end{proof}

    \section{Blocking rational directions}
    \label{sec:rational-blocking}

    Orthogonal projection onto the perpendicular plane turns blocking in a
    fixed rational direction into a two-dimensional lattice-covering problem.
    To extend the result to arbitrary directions, we will approximate each
    direction by a rational one. The resulting change of direction must
    preserve a collision, so we need a positive blocking margin. We therefore
    prove that blocking persists after every scatterer radius is reduced by
    the same positive amount.

    A \emph{uniform blocking margin} for a specified family of lines is a
    positive amount by which every scatterer radius can be reduced while
    every line in that family still enters the interior of a reduced ball.

    Throughout this section, put
    \[
        \BlockingMargin=\frac1{250}.
    \]
    For each scatterer in \(\mathcal B\) with center \(\vect{c}\)
    and radius \(r\), define the corresponding \emph{witness ball} to be the
    closed ball with center \(\vect{c}\) and radius \(r-\BlockingMargin\).
    We prove that the witness balls block every rational direction.
    The value \(\BlockingMargin=1/250\) provides a convenient
    common margin verified by exact arithmetic. No optimality is claimed.

    \subsection{Primitive vectors and the projection criterion}
    \label{subsec:primitive-periods}

    Orthogonal projection sends parallel lines to points and the interiors
    of FCC witness balls to open disks. We introduce the notation needed to
    express blocking as a covering problem; the supporting lattice facts
    are proved in Appendix~\ref{app:lattice-facts}.

    The lattices \(2\Z^3\), \(\FCC\), and \(\Z^3\) have the same rational
    directions. Indeed,
    \[
        2\Z^3\subset\FCC\subset\Z^3,
    \]
    and doubling any integer vector gives a vector of \(2\Z^3\).

    A nonzero vector \(\vect{w}\) in a lattice \(\latt{L}\) is
    \emph{primitive} if \(\latt{L}\cap\R\vect{w}=\Z\vect{w}\), or
    equivalently if its integer coordinates in a lattice basis have greatest
    common divisor one. Such a vector extends to a lattice basis
    (Appendix~\ref{app:lattice-bases}).

    Fix a nonzero vector \(\vect{v}\) representing a rational direction.
    Choose and fix \(\vect{a}\in\Z^3\) parallel to
    \(\vect{v}\) with \(\gcd(a_1,a_2,a_3)=1\). The vector \(\vect{a}\)
    is a \emph{primitive integer representative} of this direction.
    The only primitive integer representatives of this direction are
    \(\vect{a}\) and \(-\vect{a}\).

    The vector \(\vect{a}\) is primitive in \(\Z^3\). For the FCC core \(\FCC\),
    define the \emph{primitive FCC vector} \(\vect{w}(\vect{a})\) as
    the shortest positive integer multiple of \(\vect{a}\) belonging
    to \(\FCC\). Since \(\FCC\) consists of integer vectors with even
    coordinate sum,
    \begin{equation}
        \vect{w}(\vect{a})=
        \begin{cases}
            \vect{a},&a_1+a_2+a_3\text{ is even},\\
            2\vect{a},&a_1+a_2+a_3\text{ is odd}.
        \end{cases}
        \label{eq:fcc-period}
    \end{equation}
    Indeed, when \(a_1+a_2+a_3\) is odd, \(\vect{a}\notin\FCC\), and the
    smallest positive integer \(k\) satisfying \(k\vect{a}\in\FCC\)
    is \(k=2\).

    For a nonzero vector \(\vect{w}\in\R^3\), put
    \(\hat{\vect{w}}=\vect{w}/\norm{\vect{w}}\) and define the perpendicular
    plane
    \[
        \orthcomp{w}
        =\{\vect{x}\in\R^3:\vect{x}\dotprod\vect{w}=0\},
    \]
    and let \(\orthproj{w}:\R^3\to\orthcomp{w}\) be orthogonal
    projection, given by
    \begin{equation}
        \orthproj{w}(\vect{x})
        =\vect{x}-(\vect{x}\dotprod\hat{\vect{w}})\hat{\vect{w}},
        \qquad \vect{x}\in\R^3.
        \label{eq:orthogonal-projection}
    \end{equation}
    The projected inner-product identity
    \eqref{eq:projected-inner-product}, derived in
    Appendix~\ref{app:lattice-bases}, will be used for the exact calculations.

    Now fix a primitive vector \(\vect{w}\in\FCC\).
    Define the projected FCC lattice by
    \[
        \projlat{w}=\orthproj{w}(\FCC).
    \]
    Choose \(\vect{b}_1,\vect{b}_2\in\FCC\) so that
    \(\vect{w},\vect{b}_1,\vect{b}_2\) form an FCC basis. The vectors
    \(\orthproj{w}(\vect{b}_1)\) and
    \(\orthproj{w}(\vect{b}_2)\) form a basis of \(\projlat{w}\).
    Thus \(\projlat{w}\) is a two-dimensional lattice in
    \(\orthcomp{w}\).
    For every nonzero integer \(k\),
    \((k\vect{w})^\perp=\orthcomp{w}\), so the projection depends only on
    the direction of \(\vect{w}\). Primitivity of \(\vect{w}\)
    is needed for the covolume formula in
    Lemma~\ref{lem:projected-covolume}.

    For a planar lattice \(\Gamma\) spanning a plane \(P\), define the
    \emph{covering radius}
    \[
        \mu(\Gamma)=\sup_{\vect{x}\in P}\dist(\vect{x},\Gamma).
    \]
    For \(r>0\), the closed radius-\(r\) disks centered at \(\Gamma\)
    cover \(P\) exactly when \(\mu(\Gamma)\le r\). Blocking requires
    entry into scatterer interiors. Accordingly,
    Lemma~\ref{lem:projection-criterion} uses open disks and the
    strict bound \(\mu(\projlat{w})<r\) to establish blocking by
    the radius-\(r\) FCC balls.

    \begin{lemma}[Projection criterion]
        \label{lem:projection-criterion}
        Let \(\vect{w}\in\FCC\) be primitive in \(\FCC\), and let \(r>0\).
        The radius-\(r\) FCC balls block every line parallel to \(\vect{w}\) if and
        only if the open radius-\(r\) disks centered at \(\projlat{w}\) cover
        \(\orthcomp{w}\). In particular, the balls block every such line whenever
        \(\mu(\projlat{w})<r\).
    \end{lemma}

    \begin{proof}
        Let \(\ell\parallel\vect{w}\), and write
        \(\vect{x}=\orthproj{w}(\ell)\) for its projected point.
        For every FCC center \(\vect{c}\),
        \[
            \dist(\vect{c},\ell)
            =\norm{\orthproj{w}(\vect{c})-\vect{x}}.
        \]
        Thus \(\ell\) enters the open radius-\(r\) ball centered at
        \(\vect{c}\) exactly when \(\vect{x}\) belongs to the corresponding
        open disk.
    \end{proof}

    \subsection{The general projected-lattice estimate}
    \label{subsec:projected-estimate}

    To reduce the proof to finitely many directions, we combine two geometric
    bounds: the projected covolume produces a short planar vector, while the
    dual lattice controls the distance to a parallel
    plane containing lattice points. The planes considered are parallel to
    a plane spanned by two independent lattice vectors. Together
    they bound the distance from any line in the chosen direction to a lattice
    center. The estimate is valid for a general full-rank lattice; its FCC
    specialization gives \(\mu(\projlat{w})<7/10-\BlockingMargin\) whenever
    \(\norm{\vect{w}}\geq4\).

    For a full-rank lattice \(\latt{L}\subset\R^n\), write
    \(\covol(\latt{L})\) for the volume of a fundamental
    parallelepiped, and \(\lambda_1(\latt{L})\)
    for the length of its shortest nonzero vector. Its dual lattice is
    \[
        \latt{L}^*=\{\vect{y}\in\R^n:\vect{y}\dotprod\vect{x}\in\Z
        \text{ for every }\vect{x}\in\latt{L}\}.
    \]
    The two bounds used below are proved in Appendix~\ref{app:lattice-facts}.

    \begin{lemma}[Projected lattice-covering estimate]
        \label{lem:projected-covering-estimate}
        Let \(\latt{L}\subset\R^3\) be a full-rank lattice and let
        \(\vect{w}\in\latt{L}\) be primitive. Then
        \begin{equation}
            \mu\bigl(\orthproj{w}(\latt{L})\bigr)^2
            \leq
            \frac{1}{4\lambda_1(\latt{L}^*)^2}
            +\frac{\covol(\latt{L})}{2\sqrt3\,\norm{\vect{w}}}.
            \label{eq:abstract-covering-estimate}
        \end{equation}
    \end{lemma}

    \begin{proof}
        By Lemmas~\ref{lem:projected-covolume}
        and~\ref{lem:planar-short-vector}, the projected lattice contains a
        nonzero vector \(\planarvec{\xi}\) satisfying
        \[
            \norm{\planarvec{\xi}}^2
            \leq\frac{2\covol(\latt{L})}
            {\sqrt3\,\norm{\vect{w}}}.
        \]
        Choose \(\vect{b}\in\latt{L}\) with
        \(\orthproj{w}(\vect{b})=\planarvec{\xi}\), and set
        \(V=\Span(\vect{w},\vect{b})\).

        To bound the covering radius of \(\orthproj{w}(\latt{L})\),
        we must show that every point of \(\orthcomp{w}\) lies close
        to a projected lattice point. Each point of \(\orthcomp{w}\)
        represents a line \(\ell\parallel\vect{w}\), and each projected
        lattice point represents a parallel line containing lattice points.
        The distance between the corresponding points equals the distance
        between the parallel lines. We therefore bound the distance from
        an arbitrary line \(\ell\parallel\vect{w}\) to a lattice line
        by two perpendicular translations (Figure~\ref{fig:two-moves}).

        \begin{figure}[htbp]
            \centering
            \includegraphics[width=\linewidth]{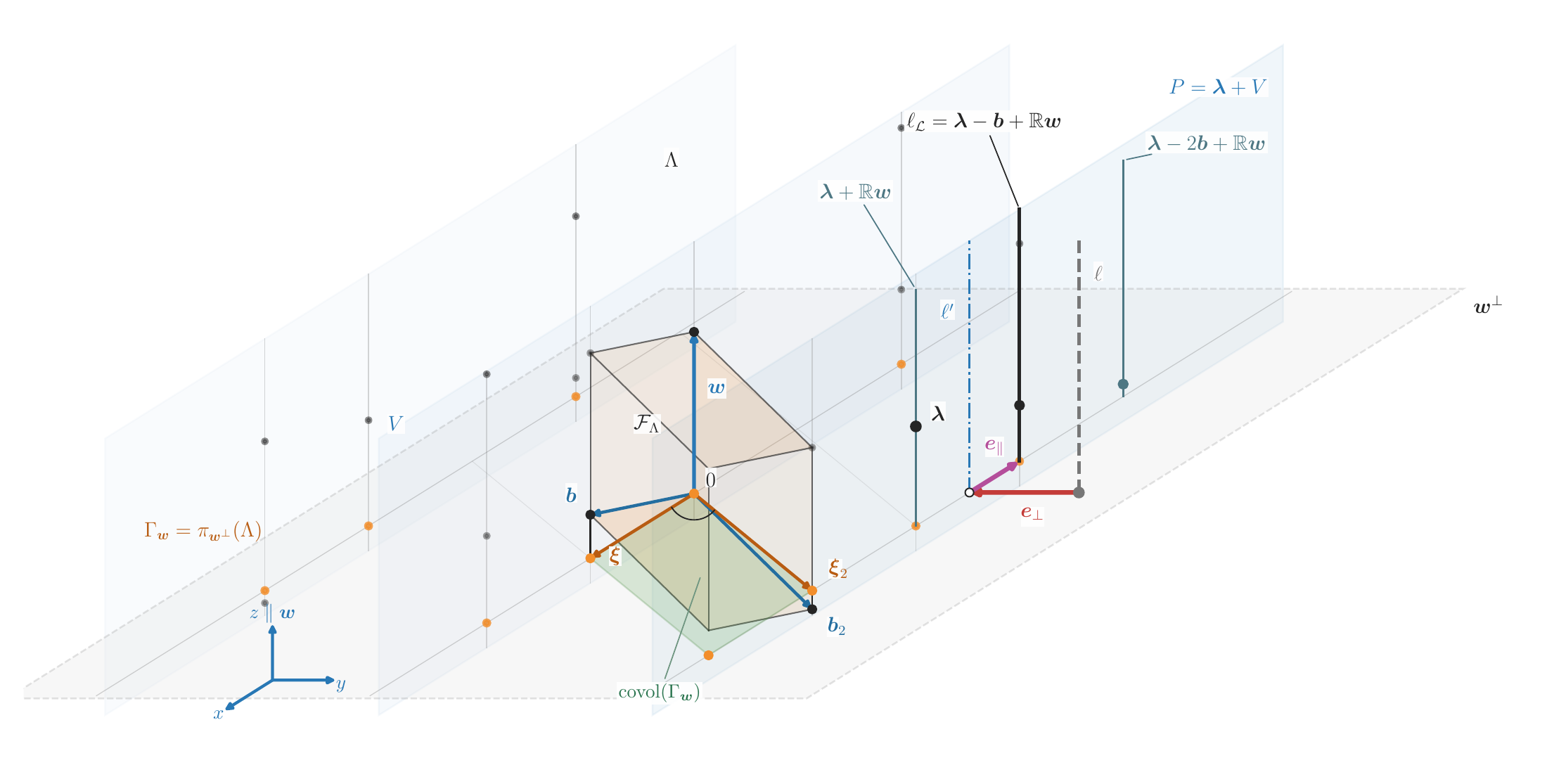}
            \caption{Two perpendicular translations carry \(\ell\) to a lattice
                line in the FCC example \(\latt{L}=\FCC\).
                The vectors \(\vect{w},\vect{b},\vect{b}_2\) span the
                fundamental parallelepiped \(\funddom{\FCC}\), and
                \(\planarvec{\xi}_2=\orthproj{w}(\vect{b}_2)\).}
            \label{fig:two-moves}
        \end{figure}

        Given a line \(\ell\parallel\vect{w}\), let \(Q\) be the plane
        containing \(\ell\) and parallel to \(V\).
        Lemma~\ref{lem:lattice-plane-spacing} provides a plane
        \(P\parallel V\) containing lattice points, with
        \[
            \dist(Q,P)\le\frac{1}{2\lambda_1(\latt{L}^*)}.
        \]
        The first move translates \(\ell\) perpendicularly into \(P\).
        Denote the translation vector by \(\vect{e}_{\perp}\) and the
        translated line by \(\ell'=\ell+\vect{e}_{\perp}\). Then
        \[
            \vect{e}_{\perp}\perp V,
            \qquad
            \norm{\vect{e}_{\perp}}
            \le\frac{1}{2\lambda_1(\latt{L}^*)}.
        \]

        Choose \(\vect{\lambda}\in\latt{L}\cap P\), so
        \(P=\vect{\lambda}+V\). The plane \(P\) contains the parallel
        lattice lines
        \[
            \vect{\lambda}+n\vect{b}+\R\vect{w},
            \qquad n\in\Z.
        \]
        Consecutive lattice lines are separated by
        \(\norm{\orthproj{w}(\vect{b})}=\norm{\planarvec{\xi}}\).
        The second move translates \(\ell'\) within \(P\), perpendicularly
        to \(\vect{w}\), to a nearest line
        \(\ell_{\latt{L}}\) in the displayed family
        \(\{\vect{\lambda}+n\vect{b}+\R\vect{w}:n\in\Z\}\).
        Writing \(\vect{e}_{\parallel}\) for the second translation vector,
        we obtain
        \[
            \ell_{\latt{L}}=\ell'+\vect{e}_{\parallel},
            \qquad
            \vect{e}_{\parallel}\in V,
            \qquad
            \norm{\vect{e}_{\parallel}}
            \le\frac{\norm{\planarvec{\xi}}}{2}.
        \]

        The two moves are orthogonal because
        \(\vect{e}_{\perp}\perp V\) and \(\vect{e}_{\parallel}\in V\).
        Both moves are perpendicular to \(\vect{w}\), so
        \[
            \dist(\ell,\ell_{\latt{L}})^2
            =\norm{\vect{e}_{\perp}}^2+\norm{\vect{e}_{\parallel}}^2
            \le\frac{1}{4\lambda_1(\latt{L}^*)^2}
            +\frac{\norm{\planarvec{\xi}}^2}{4}.
        \]
        Since \(\ell_{\latt{L}}\) contains a lattice point,
        \[
            \dist\bigl(\orthproj{w}(\ell),\orthproj{w}(\latt{L})\bigr)
            \le \dist(\ell,\ell_{\latt{L}}).
        \]
        As \(\ell\) ranges over all lines parallel to \(\vect{w}\),
        \(\orthproj{w}(\ell)\) ranges over \(\orthcomp{w}\).
        Taking the supremum and using the preceding bound on
        \(\norm{\planarvec{\xi}}^2\) proves
        \eqref{eq:abstract-covering-estimate}.
    \end{proof}

    \begin{remark}
        In~\eqref{eq:abstract-covering-estimate}, the second term tends to zero as
        \(\norm{\vect{w}}\to\infty\), while the first depends only on the spacing
        of parallel planes containing lattice points
        and parallel to \(\Span(\vect{u},\vect{v})\), where
        \(\vect{u},\vect{v}\in\latt{L}\) are linearly independent.
        Therefore a ball radius above the plane-spacing
        threshold leaves only finitely many directions with short primitive
        lattice vectors to check.
    \end{remark}

    \subsection{Specialization to the FCC core}

    The FCC lattice has covolume \(2\), since
    \[
        \Z^3=\FCC\,\sqcup\,\bigl((1,0,0)+\FCC\bigr).
    \]
    Its dual lattice and minimum are
    \[
        \dualFCC=\Z^3\cup\bigl((1/2,1/2,1/2)+\Z^3\bigr),
        \qquad \lambda_1(\dualFCC)=\frac{\sqrt3}{2},
    \]
    as proved in Lemma~\ref{lem:fcc-dual} of
    Appendix~\ref{app:lattice-facts}.

    The packing and plane-spacing thresholds explain the choice of radius \(7/10\) for the balls centered at \(\FCC\).
    The shortest nonzero FCC vectors have length \(\sqrt2\), and hence equal
    lattice-centered balls form a strictly separated packing provided
    \(R<\lambda_1(\FCC)/2=1/\sqrt2\). On the other hand,
    Lemma~\ref{lem:lattice-plane-spacing} shows that every plane parallel
    to a plane spanned by two independent FCC vectors meets the interior
    of an FCC-centered \(R\)-ball when
    \(R>1/(2\lambda_1(\dualFCC))=1/\sqrt3\). Since
    \[
        \frac1{\sqrt3}<\frac7{10}-\BlockingMargin<\frac7{10}<\frac1{\sqrt2},
    \]
    even the FCC witness radius remains above the plane-spacing
    threshold, while the actual radius remains below the packing threshold.

    Substituting the FCC data into
    Lemma~\ref{lem:projected-covering-estimate} gives a blocking criterion
    in terms of the length of the primitive FCC vector \(\vect{w}\).

    \begin{lemma}[FCC line estimate]
        \label{lem:fcc-lines}
        Let \(\vect{w}\in\FCC\) be primitive in \(\FCC\). If
        \(\norm{\vect{w}}\ge4\), every line parallel to \(\vect{w}\) enters an FCC witness ball of radius \(7/10-\BlockingMargin\).
    \end{lemma}

    \begin{proof}
        Substitute \(\covol(\FCC)=2\) and
        \(\lambda_1(\dualFCC)=\sqrt3/2\) into the abstract
        estimate~\eqref{eq:abstract-covering-estimate} to obtain
        \[
            \mu(\projlat{w})^2
            \leq\frac13+\frac1{\sqrt3\,\norm{\vect{w}}}.
        \]
        If \(\norm{\vect{w}}\ge4\), a direct numerical comparison gives
        \[
            \mu(\projlat{w})^2
            \le\frac13+\frac1{4\sqrt3}
            <\left(\frac{174}{250}\right)^2
            =\left(\frac7{10}-\BlockingMargin\right)^2.
        \]
        Thus \(\mu(\projlat{w})<7/10-\BlockingMargin\).
        By Lemma~\ref{lem:projection-criterion}, the FCC witness balls
        block every line parallel to \(\vect{w}\).
    \end{proof}

    \subsection{Exceptional directions and the finite certificate}

    Lemma~\ref{lem:fcc-lines} proves that the FCC witness balls
    block every direction whose primitive FCC vector \(\vect{w}\)
    satisfies \(\norm{\vect{w}}\ge4\). It remains to prove blocking
    for directions whose primitive FCC vectors have length less than \(4\).

    The \emph{type} of a primitive FCC vector \(\vect{w}\) is the triple
    \((x,y,z)\) obtained by arranging the absolute values of the coordinates
    of \(\vect{w}\) in the order \(x\ge y\ge z\ge0\).
    Vectors of the same type differ only by coordinate permutations and
    sign changes. Signed coordinate permutations preserve \(\FCC\) and
    Euclidean distances, so vectors of the same type have isometric projected
    lattices and equal covering radii. It therefore suffices to compute one
    covering radius for each type.
    The enumeration in Lemma~\ref{lem:directions} gives 61 unoriented
    directions, grouped into the six types in
    Table~\ref{tab:short-fcc-radii}.

    For each representative vector \(\vect{w}\) in
    Table~\ref{tab:short-fcc-radii}, the listed vectors
    \(\vect{w},\vect{b}_1,\vect{b}_2\) form an FCC basis. Thus
    \(\planarvec{\xi}_i=\orthproj{w}(\vect{b}_i)\), \(i=1,2\), form a
    basis of \(\projlat{w}\). Write its Gram triple as
    \[
        (g_{11},g_{12},g_{22})
        =\bigl(\norm{\planarvec{\xi}_1}^2,
        \planarvec{\xi}_1\dotprod\planarvec{\xi}_2,
        \norm{\planarvec{\xi}_2}^2\bigr).
    \]
    The listed bases satisfy \(g_{11}\le g_{22}\) and
    \(-g_{11}/2\le g_{12}\le0\). Under these inequalities,
    Lemma~\ref{lem:planar-covering} in Appendix~\ref{app:planar-covering}
    gives
    \[
        \mu(\projlat{w})^2
        =\frac{g_{11}g_{22}(g_{11}+g_{22}+2g_{12})}
        {4(g_{11}g_{22}-g_{12}^2)}.
    \]
    Geometrically, this is the squared circumradius of either triangle
    obtained by dividing the basis parallelogram along the diagonal from
    \(0\) to \(\planarvec{\xi}_1+\planarvec{\xi}_2\).
    The appendix proves the formula and supplies Gauss reduction for
    choosing such bases. Here the bases are explicit, so it remains to
    verify their Gram entries and compare the resulting covering radii
    with \(7/10-\BlockingMargin\). This comparison settles 52 of the
    61 short directions.

    \begin{table}[htbp]
        \centering
        \scriptsize
        \setlength{\tabcolsep}{3.2pt}
        \begin{tabular}{c c c c c c c}
            \toprule
            type \(\vect{w}\) & \(\vect{b}_1\) & \(\vect{b}_2\) &
            Gram triple & directions & \(\mu(\projlat{w})^2\) & conclusion \\
            \midrule
            \((1,1,0)\) & \((-1,0,-1)\) & \((-1,0,1)\) &
            \((3/2,-1/2,3/2)\) & 6 & \(9/16\) & correction \\
            \((2,0,0)\) & \((1,0,1)\) & \((1,-1,0)\) &
            \((1,0,1)\) & 3 & \(1/2\) & correction \\
            \((2,1,1)\) & \((1,0,1)\) & \((0,-1,-1)\) &
            \((1/2,0,4/3)\) & 12 & \(11/24\) & core \\
            \((2,2,2)\) & \((0,1,1)\) & \((1,0,1)\) &
            \((2/3,-1/3,2/3)\) & 4 & \(2/9\) & core \\
            \((3,1,0)\) & \((-1,-1,0)\) & \((-1,0,1)\) &
            \((2/5,-1/5,11/10)\) & 12 & \(121/400\) & core \\
            \((3,2,1)\) & \((1,1,0)\) & \((0,1,1)\) &
            \((3/14,-1/14,19/14)\) & 24 & \(285/784\) & core \\
            \bottomrule
        \end{tabular}
        \caption{Exact projected FCC covering radii for the six short-vector types.
        The directions column counts unoriented directions of each type.
        ``Core'' means that the FCC witness balls suffice.
        ``Correction'' refers to the joint covering
        certificate for the FCC core and the three smaller-ball families.}
        \label{tab:short-fcc-radii}
    \end{table}

    \begin{lemma}[Short directions]
        \label{lem:directions}
        Exactly 61 unoriented rational directions have
        \(\norm{\vect{w}(\vect{a})}<4\). The FCC witness balls already
        block 52 of them. The remaining nine directions have primitive FCC
        vectors of type \((1,1,0)\) or \((2,0,0)\).
    \end{lemma}

    \begin{proof}
        Write the type as \((x,y,z)\), with \(x\ge y\ge z\ge0\). The bound
        \(x^2+y^2+z^2<16\) gives \(x\le3\), and the FCC condition makes
        \(x+y+z\) even. If \(x=1\), only \((1,1,0)\) occurs. If \(x=2\), the
        possibilities are \((2,0,0)\), \((2,1,1)\), \((2,2,0)\), and
        \((2,2,2)\); the vector \((2,2,0)=2(1,1,0)\) is not primitive in
        \(\FCC\). If \(x=3\), the norm and parity conditions leave
        \((3,1,0)\) and \((3,2,1)\). Thus the six primitive types are
        \[
            (1,1,0),\quad (2,0,0),\quad (2,1,1),\quad
            (2,2,2),\quad (3,1,0),\quad (3,2,1).
        \]
        For each row of Table~\ref{tab:short-fcc-radii}, the FCC vectors
        \(\vect{w},\vect{b}_1,\vect{b}_2\) have determinant
        \(\pm2=\covol(\FCC)\) and therefore form a basis of \(\FCC\).
        The Gram entries for the projected basis vectors
        \(\planarvec{\xi}_i=\orthproj{w}(\vect{b}_i)\) satisfy
        \(g_{11}\le g_{22}\) and \(-g_{11}/2\le g_{12}\le0\), as in
        \eqref{eq:gauss-reduced-basis}. Applying
        Lemma~\ref{lem:planar-covering} with \((g_{11},g_{12},g_{22})\) gives the exact
        radii in Table~\ref{tab:short-fcc-radii}.

        For each type, count the distinct vectors obtained by permuting
        coordinates and independently choosing signs for the nonzero coordinates.
        Each unoriented direction corresponds to a pair
        \(\{\vect{u},-\vect{u}\}\), so the number of directions is half
        the number of vectors. In the order of Table~\ref{tab:short-fcc-radii},
        the direction counts are
        \[
            3\cdot4/2=6,\quad 3\cdot2/2=3,\quad
            3\cdot8/2=12,\quad 8/2=4,\quad
            6\cdot4/2=12,\quad 6\cdot8/2=24.
        \]
        Their sum is \(61\). The largest squared covering radius among the last
        four types is \(11/24\), and
        \[
            \frac{11}{24}<\left(\frac7{10}-\BlockingMargin\right)^2.
        \]
        Hence the FCC witness balls block every line in those 52 directions. Only the
        first two types, totaling nine directions, require correction balls.
    \end{proof}

    Signed coordinate permutations preserve the FCC core, which allowed us
    to compute one covering radius per type. We do not assume these
    symmetries for the full configuration, so the certificate checks each
    of the nine exceptional directions separately.

    For each of the nine directions remaining in Lemma~\ref{lem:directions},
    fix a primitive FCC vector \(\vect{w}\). The core estimates used
    \(\orthproj{w}(\FCC)\). To include all four ball families, we now work
    modulo the specified period lattice \(2\Z^3\), so the relevant
    \emph{line space} is
    \[
        \mathcal S_{\vect{w}}
        =\orthcomp{w}/\orthproj{w}(2\Z^3).
    \]
    The points of \(\mathcal S_{\vect{w}}\) represent lines parallel
    to \(\vect{w}\), with lines identified under translations in
    \(2\Z^3\). Since \(2\Z^3\subset\FCC\),
    \(\orthproj{w}(2\Z^3)\subset\projlat{w}\).

    A \emph{projected witness disk} is the open disk obtained by orthogonally
    projecting the interior of a witness ball onto \(\orthcomp{w}\).
    In line space, the term denotes the image of this disk under the quotient
    map. We certify blocking by covering the compact line space
    \(\mathcal S_{\vect{w}}\) with projected witness disks.

    The certificate subdivides a parameter square representing
    \(\mathcal S_{\vect{w}}\) into dyadic cells. Each parameter cell represents
    a parallelogram in the perpendicular plane. Suppose its center is
    \(\planarvec{p}\), and every point of the parallelogram lies within
    distance \(\rho\) of \(\planarvec{p}\). A projected witness disk
    centered at \(\planarvec{q}\), of radius \(r-\BlockingMargin\), contains
    the entire cell whenever
    \[
        \norm{\planarvec{p}-\planarvec{q}}+\rho<r-\BlockingMargin.
    \]
    Appendix~\ref{app:certificate} converts this sufficient condition into
    strict integer inequalities. A cell is accepted only after verifying
    this condition for one witness disk; otherwise it is divided into four
    smaller cells. If subdivision terminates with every terminal cell
    accepted, the disks cover \(\mathcal S_{\vect{w}}\), so every line
    parallel to \(\vect{w}\) is blocked by a witness ball.

    \begin{lemma}[Exact line-space certificate]
        \label{lem:certificate}
        For every primitive FCC vector \(\vect{w}\) in one of the nine
        directions identified in Lemma~\ref{lem:directions},
        the open disks obtained by projecting the interiors of the 38 witness balls cover the line space \(\mathcal S_{\vect{w}}\).
    \end{lemma}

    \begin{proof}
        Using the cell-acceptance test derived in
        Appendix~\ref{app:certificate}, the algorithm in
        Appendix~\ref{alg:line-certificate} tests a total of 9213 cells and
        accepts 6912 across the nine directions. The maximum subdivision
        depth over all nine runs is \(11\).
        Each accepted cell lies in a projected witness disk by that test.
        Subdivision preserves coverage of the initial parameter square, and
        termination leaves only accepted terminal cells. Hence the disks cover
        \(\mathcal S_{\vect{w}}\) for each of the nine directions.
    \end{proof}

    \begin{lemma}[Uniform margin in rational directions]
        \label{lem:rational-margin}
        Every line with rational direction enters the interior of one of the
        witness balls obtained by reducing the radius of each ball in
        \(\mathcal B\) by \(\BlockingMargin=1/250\).
    \end{lemma}

    \begin{proof}
        Let \(\vect{a}\in\Z^3\) be the primitive representative of the direction.
        If \(\norm{\vect{w}(\vect{a})}\geq4\), Lemma~\ref{lem:fcc-lines} applies.
        If the norm is less than \(4\), Lemma~\ref{lem:directions} handles 52 of
        these 61 directions using the FCC witness balls, and
        Lemma~\ref{lem:certificate} handles the remaining nine.
    \end{proof}

    Lemma~\ref{lem:rational-margin} establishes blocking for every rational
    direction with the common margin \(\BlockingMargin\).
    Section~\ref{sec:quantitative} turns the margin \(\BlockingMargin\)
    into a uniform bound on the lengths of collision-free segments in all
    directions, completing the proof of finite horizon.

    \section{An explicit upper bound for the horizon}
    \label{sec:quantitative}

    The uniform margin in Lemma~\ref{lem:rational-margin} lets us replace an
    arbitrary direction by a nearby rational one while retaining a collision.
    The approximation must control both the angular error and the period
    \(2\norm{\vect{a}}\) along a line with primitive integer representative
    \(\vect{a}\). Along such a line, a witness-ball intersection occurs within
    distance \(\norm{\vect{a}}\) of any prescribed point.
    Transferring this intersection
    to a line through that point making angle \(\theta\) introduces a
    displacement of at most \(\norm{\vect{a}}\tan\theta\), as shown in
    Lemma~\ref{lem:angular-transfer}. We therefore seek an approximation for
    which this displacement is smaller than \(\BlockingMargin\), while
    \(\norm{\vect{a}}\) remains uniformly bounded. The rational approximation
    step uses the following simultaneous form of Dirichlet's theorem.

    \begin{lemma}[Simultaneous Dirichlet approximation]
        \label{lem:dirichlet}
        For real numbers \(x,y\) and an integer \(N\ge1\), there are integers
        \(q,p,r\) such that
        \[
            1\le q\le N^2,
            \qquad
            \abs{qx-p}\le\frac1N,
            \qquad
            \abs{qy-r}\le\frac1N.
        \]
    \end{lemma}

    \begin{proof}
        Consider the \(N^2+1\) points
        \[
            (\{kx\},\{ky\}),\qquad 0\le k\le N^2,
        \]
        in the unit square and divide the square into \(N^2\) squares of side
        \(1/N\). Two points, with indices \(k<\ell\), lie in the same small
        square. Taking \(q=\ell-k\), the two coordinate differences are within
        \(1/N\) of integers \(p\) and \(r\), respectively. This gives all three
        stated inequalities.
    \end{proof}

    \subsection{Upper bound by rational approximation}

    Along a rational line, periodicity bounds the distance from
    any prescribed point to a witness-ball intersection. The following lemma
    uses the distance bound and the radius reduction to retain a collision
    when the direction is changed.

    \begin{lemma}[Collision transfer from a nearby rational direction]
        \label{lem:angular-transfer}
        Let \(\vect{a}\in\Z^3\) be primitive, let
        \(\vect{v}\in S^2\), and let \(0\le\theta<\pi/2\) be the angle
        between the lines \(\R\vect{a}\) and \(\R\vect{v}\).
        Suppose every line parallel to \(\vect{a}\) meets a
        concentric witness ball obtained
        by decreasing each scatterer radius by \(\eta>0\). If
        \[
            \norm{\vect{a}}\tan\theta<\eta,
        \]
        then every line parallel to \(\vect{v}\)
        enters the interior of an actual scatterer within
        distance at most \(\norm{\vect{a}}/\cos\theta\) of any prescribed point on
        the line.
    \end{lemma}

    \begin{proof}
        Fix an arbitrary prescribed point on a line parallel to
        \(\vect{v}\), and translate coordinates so that the prescribed point
        is the origin. The period lattice remains \(2\Z^3\), and the blocking
        hypothesis is unchanged.
        Choose the sign of \(\vect{a}\) so that the unit vector
        \(\hat{\vect{a}}=\vect{a}/\norm{\vect{a}}\) satisfies
        \(\hat{\vect{a}}\dotprod\vect{v}=\cos\theta>0\), and put
        \(\ell_a=\R\hat{\vect{a}}\) and \(\ell_v=\R\vect{v}\).
        By hypothesis, \(\ell_a\) meets a witness ball. Translation by
        \(2\vect{a}\) preserves both \(\ell_a\) and the periodic family of
        witness balls. Here \(2\vect{a}\) is a period
        of the full configuration; \(\vect{w}(\vect{a})\) was the primitive
        FCC vector used for the core. The intersections therefore repeat every
        \(2\norm{\vect{a}}\) along \(\ell_a\), so we may choose an intersection
        point \(\vect{x}_a\in\ell_a\) within half a period of the origin.
        For some \(t\in\R\) and a scatterer
        of radius \(r\) with translated center \(\vect{c}\), we have
        \[
            \vect{x}_a=t\hat{\vect{a}},\qquad
            \abs{t}\leq\norm{\vect{a}},\qquad
            \norm{\vect{x}_a-\vect{c}}\le r-\eta.
        \]

        Choose \(\vect{x}_v\in\ell_v\) so that
        \(\vect{x}_v-\vect{x}_a\perp\hat{\vect{a}}\)
        (Figure~\ref{fig:collision-transfer}).

        \begin{figure}[htbp]
            \centering
            \includegraphics[width=\linewidth]{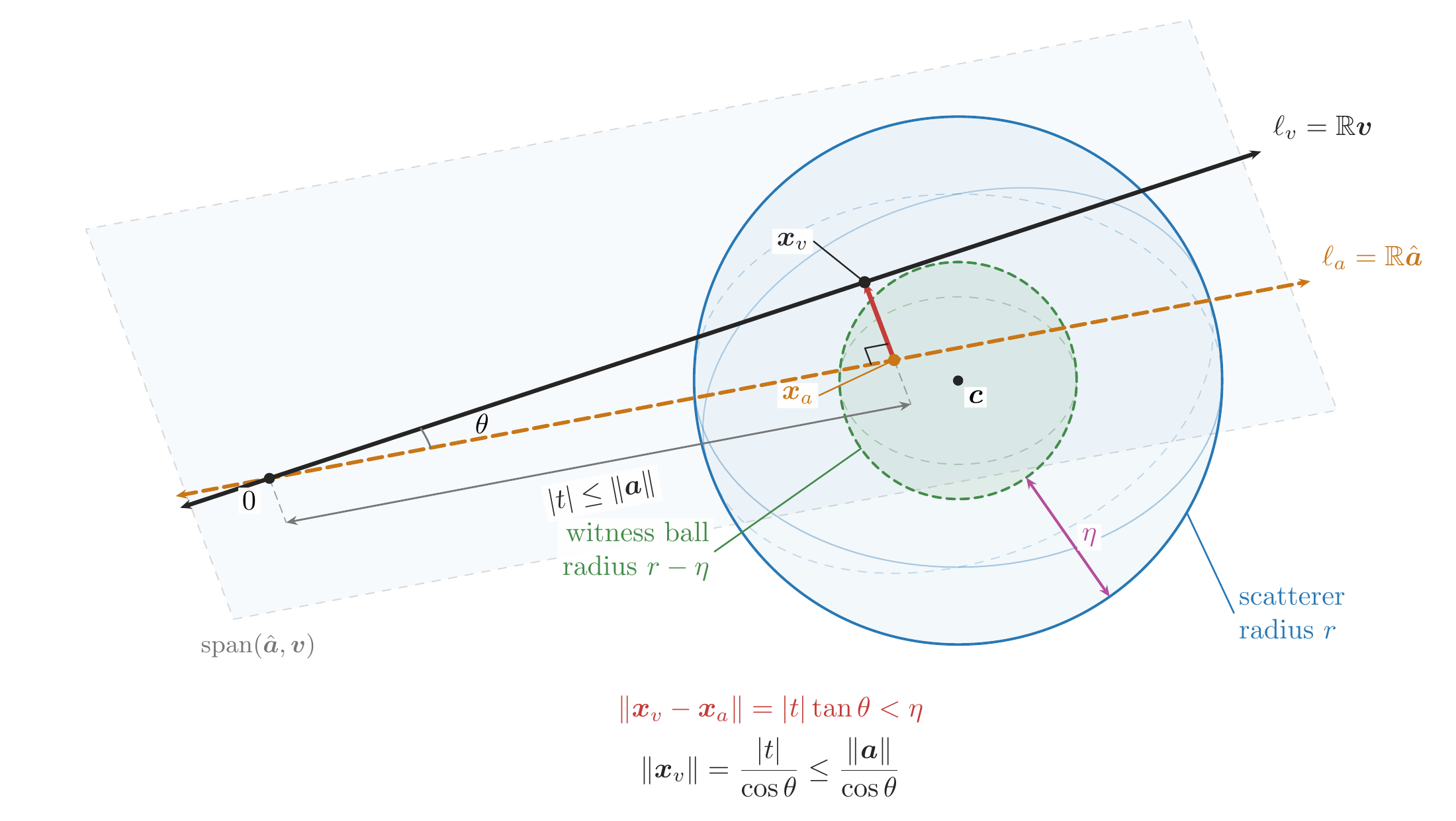}
            \caption{The point \(\vect{x}_a\) lies in the witness ball. The
                displacement to \(\vect{x}_v\) is less than \(\eta\), so
                \(\vect{x}_v\) lies in the interior of the scatterer.}
            \label{fig:collision-transfer}
        \end{figure}

        The perpendicularity condition gives
        \(\vect{x}_v\dotprod\hat{\vect{a}}
        =\vect{x}_a\dotprod\hat{\vect{a}}=t\).
        Since \(\vect{v}\dotprod\hat{\vect{a}}=\cos\theta\),
        \[
            \vect{x}_v=\frac{t}{\cos\theta}\vect{v}.
        \]
        The vectors \(\vect{x}_a\) and \(\vect{x}_v-\vect{x}_a\) are
        perpendicular, so the Pythagorean theorem gives
        \begin{align*}
            \norm{\vect{x}_v-\vect{x}_a}^2
            &=\frac{t^2}{\cos^2\theta}-t^2
            =t^2\tan^2\theta,\\
            \norm{\vect{x}_v-\vect{x}_a}
            &=\abs{t}\tan\theta
            \leq\norm{\vect{a}}\tan\theta<\eta.
        \end{align*}
        Consequently,
        \[
            \norm{\vect{x}_v-\vect{c}}
            \leq\norm{\vect{x}_v-\vect{x}_a}
            +\norm{\vect{x}_a-\vect{c}}
            <\eta+(r-\eta)=r,
        \]
        so \(\vect{x}_v\in\ell_v\) lies in the actual scatterer interior.
        The distance from the prescribed point to \(\vect{x}_v\) satisfies
        \[
            \norm{\vect{x}_v}
            =\frac{\abs{t}}{\cos\theta}
            \leq\frac{\norm{\vect{a}}}{\cos\theta}.
        \]
    \end{proof}

    To apply Lemma~\ref{lem:angular-transfer} to an arbitrary unit
    direction \(\vect{v}\), we choose a primitive integer vector \(\vect{a}\)
    such that \(\norm{\vect{a}}\tan\theta<\BlockingMargin\), where
    \(\theta\) is the angle between the lines \(\R\vect{a}\) and
    \(\R\vect{v}\). The following lemma constructs \(\vect{a}\) and
    bounds \(\norm{\vect{a}}\) uniformly in \(\vect{v}\).

    \begin{lemma}[Rational approximation with controlled angle and length]
        \label{lem:controlled-rational-approximation}
        For every \(\vect{v}\in S^2\), there exists a primitive
        integer vector \(\vect{a}\) such that the angle \(\theta\) between
        \(\R\vect{a}\) and \(\R\vect{v}\) belongs to \([0,\pi/2)\) and
        \[
            \norm{\vect{a}}\tan\theta<\BlockingMargin,\qquad
            \norm{\vect{a}}<226800+\frac1{252}.
        \]
    \end{lemma}

    \begin{proof}
        Set \(\eta=\BlockingMargin=1/250\) and
        \(N=360\), and fix \(\vect{v}\in S^2\). This choice of \(N\) gives
        \(\sqrt2/N<1/252\), which makes the displacement estimate
        in~\eqref{eq:direction-weighted-angle} smaller than \(\eta\).
        We first construct an integer vector
        \(\vect{b}\) close to the ray in direction \(\vect{v}\). The
        perpendicular distance from \(\vect{b}\) to the ray will control
        the angle between \(\vect{b}\) and \(\vect{v}\).

        Since coordinate permutations preserve integer vectors, lengths,
        and angles, relabel the coordinate axes so that
        \(\abs{v_3}=\max_i\abs{v_i}\). Since \(\vect{v}\) and
        \(-\vect{v}\) determine the same lines, choose the sign of
        \(\vect{v}\) so that \(v_3>0\). To describe the direction
        \(\vect{v}\) by two bounded coordinate ratios, put
        \[
            \kappa=\frac1{v_3},\qquad \tilde{\vect{v}}=\kappa\vect{v}=(x,y,1).
        \]
        Then \(\abs{x},\abs{y}\le1\), and
        \[
            1\le \kappa=\norm{\tilde{\vect{v}}}\le\sqrt3.
        \]
        By Lemma~\ref{lem:dirichlet}, there exist integers \(p,r,q\) such that
        \[
            1\le q\le N^2,\qquad
            \abs{qx-p}\le\frac1N,\qquad
            \abs{qy-r}\le\frac1N.
        \]
        Define
        \[
            \vect{b}=(p,r,q)\in\Z^3,\qquad
            \vect{e}=\vect{b}-q\tilde{\vect{v}}=(p-qx,r-qy,0).
        \]
        Figure~\ref{fig:rational-approximation} shows the rescaling and the
        approximation by \(\vect{b}\).

        \begin{figure}[htbp]
            \centering
            \includegraphics[width=.85\linewidth]{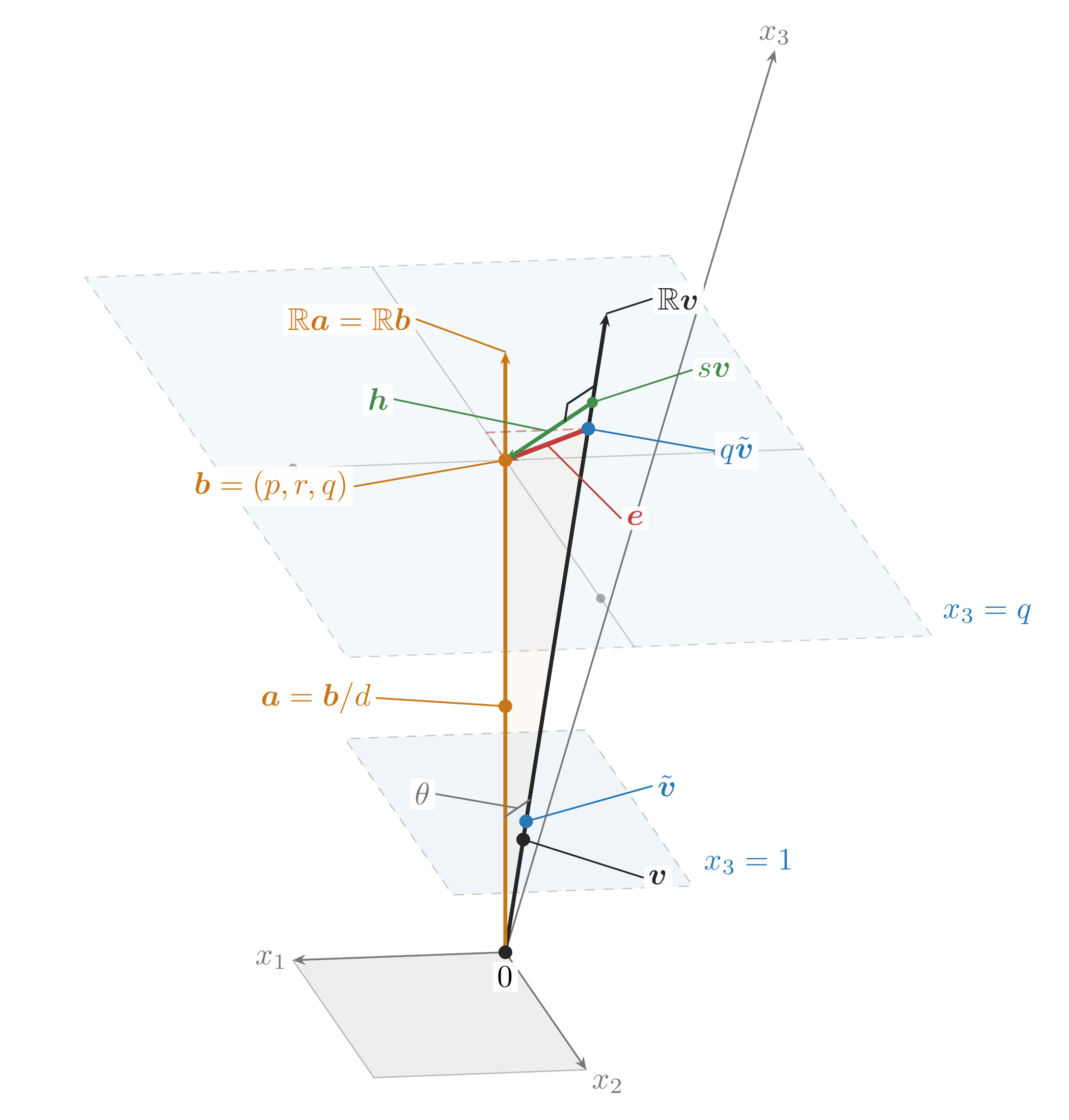}
            \caption{Rescaling \(\vect{v}\) to \(\tilde{\vect{v}}\) and
                approximating \(q\tilde{\vect{v}}\) by the integer vector
                \(\vect{b}\). The green vector \(\vect{h}\) is the
                perpendicular displacement from \(s\vect{v}\) to \(\vect{b}\).}
            \label{fig:rational-approximation}
        \end{figure}

        The component bounds therefore give
        \begin{equation}
            \norm{\vect{e}}\le\frac{\sqrt2}{N}
            =\frac{\sqrt2}{360}<\frac1{252}.
            \label{eq:direction-error-bound}
        \end{equation}
        The triangle inequality gives
        \begin{equation}
            \norm{\vect{b}}\le q\kappa+\norm{\vect{e}}.
            \label{eq:direction-vector-length}
        \end{equation}
        Let \(\theta\) be the angle between \(\vect{b}\) and \(\vect{v}\).
        Put \(s=\vect{b}\dotprod\vect{v}\), so \(s\vect{v}\) is the orthogonal
        projection of \(\vect{b}\) onto \(\R\vect{v}\), and let
        \(\vect{h}=\vect{b}-s\vect{v}\) be the perpendicular displacement.
        Since \(\vect{b}=q\kappa\vect{v}+\vect{e}\), the perpendicular
        distance \(\norm{\vect{h}}=\norm{\vect{b}}\sin\theta\) from
        \(\vect{b}\) to \(\R\vect{v}\) is at most \(\norm{\vect{e}}\).
        The forward component \(s=\norm{\vect{b}}\cos\theta\)
        equals \(q\kappa+\vect{e}\dotprod\vect{v}\), which by
        Cauchy--Schwarz is at least \(q\kappa-\norm{\vect{e}}>0\).
        Combining the sine and cosine bounds with
        \eqref{eq:direction-vector-length} gives
        \[
            \norm{\vect{b}}\tan\theta
            =\bigl(\norm{\vect{b}}\sin\theta\bigr)
              \frac{\norm{\vect{b}}}{\norm{\vect{b}}\cos\theta}
            \le \norm{\vect{e}}
                \frac{q\kappa+\norm{\vect{e}}}{q\kappa-\norm{\vect{e}}}.
        \]
        Divide the numerator and denominator of the last fraction by
        \(q\kappa\). Since \(q\kappa\ge1\), the error
        bound~\eqref{eq:direction-error-bound} gives
        \(\norm{\vect{e}}/(q\kappa)<1/252\). Hence
        \begin{equation}
            \norm{\vect{b}}\tan\theta
            <\frac1{252}\frac{1+1/252}{1-1/252}
            <\frac1{250}=\eta.
            \label{eq:direction-weighted-angle}
        \end{equation}

        Finally, set \(d=\gcd(p,r,q)\) and \(\vect{a}=\vect{b}/d\).
        Since \(q\ge1\), we have \(d\ge1\). The vector \(\vect{a}\) is
        primitive, points in the same direction as \(\vect{b}\), and
        satisfies \(\norm{\vect{a}}\le\norm{\vect{b}}\). Consequently,
        \eqref{eq:direction-weighted-angle} gives
        \(\norm{\vect{a}}\tan\theta<\eta\).
        For the length estimate, use \(q\le N^2=360^2\),
        \(\kappa\le\sqrt3<7/4\), and \eqref{eq:direction-error-bound}
        in~\eqref{eq:direction-vector-length} to obtain
        \[
            \norm{\vect{a}}\le\norm{\vect{b}}
            \le q\kappa+\norm{\vect{e}}
            <360^2\frac74+\frac1{252}
            =226800+\frac1{252}.
        \]
    \end{proof}

    We now apply Lemma~\ref{lem:angular-transfer} using the
    rational approximation from Lemma~\ref{lem:controlled-rational-approximation}
    and the blocking margin from Lemma~\ref{lem:rational-margin}.

    \begin{proposition}[Upper horizon bound]
        \label{prop:horizon-upper}
        The billiard defined by \(\mathcal B\) satisfies \(\Horizon<453601\).
    \end{proposition}

    \begin{proof}
        Set \(\eta=\BlockingMargin\). For any \(\vect{v}\in S^2\),
        choose \(\vect{a}\) using Lemma~\ref{lem:controlled-rational-approximation},
        and let \(\theta\) denote the angle between \(\R\vect{a}\) and
        \(\R\vect{v}\). Lemma~\ref{lem:rational-margin} supplies the
        witness-ball intersections required by Lemma~\ref{lem:angular-transfer}
        with radius reduction \(\eta\). Lemma~\ref{lem:angular-transfer}
        therefore ensures that every line parallel to \(\vect{v}\) enters
        a scatterer interior within distance at most
        \(\norm{\vect{a}}/\cos\theta\) of any prescribed point on the line.
        The two bounds in Lemma~\ref{lem:controlled-rational-approximation} give
        \begin{align*}
            \frac{\norm{\vect{a}}}{\cos\theta}
            &=\sqrt{\norm{\vect{a}}^2+(\norm{\vect{a}}\tan\theta)^2}\\
            &<\sqrt{\norm{\vect{a}}^2+\eta^2}<\norm{\vect{a}}+\eta\\
            &<226800+\frac1{252}+\frac1{250}=:B.
        \end{align*}
        The constant \(B\) is independent of \(\vect{v}\).

        Now consider any segment of length at least
        \(2B\). The line through the segment meets a scatterer interior less
        than \(B\) from the segment's midpoint.
        The segment therefore meets that interior and is
        not collision-free. Hence
        \[
            \Horizon\le2B
            =2\left(226800+\frac1{252}+\frac1{250}\right)<453601.
        \]

    \end{proof}

    \begin{proof}[Proof of Theorem~\ref{thm:explicit}]
        The center lists give 38 scatterers with rational centers and four
        rational radii. Proposition~\ref{prop:packing-period} gives the full
        translation lattice and strict separation. The latter makes all boundary
        components disjoint smooth spheres, so there are no corners or cusps.
        Proposition~\ref{prop:horizon-upper} gives the stated upper horizon
        bound and, in particular, finite horizon.
    \end{proof}

    To prove stability under independent perturbations, we now obtain a
    common radius reduction such that every segment
    of a fixed length still enters a reduced ball's interior. Compactness
    supplies this reduction for the original periodic configuration.

    \begin{lemma}[Uniform margin for segments]
        \label{lem:uniform-segment-margin}
        Let a periodic family of balls with finitely many
        balls modulo its period lattice have finite horizon. There are \(L>0\)
        and \(\eta>0\) such that every segment of length \(L\) meets the
        interior of a ball obtained by reducing an original ball's radius by
        \(\eta\).
    \end{lemma}

    \begin{proof}
        Choose \(L\) greater than the horizon. Every
        segment of length \(L\) then enters an original ball. Modulo the period
        lattice, such segments are parametrized by their initial points and
        directions in the compact space \(\mathbb T\times S^2\), where
        \(\mathbb T\) is the corresponding quotient torus.

        For each ball orbit of radius \(r\), \(t\in[0,L]\), and
        \(0<\delta<r\), let \(U\) be the set of initial points and directions
        for which the point at distance \(t\) along the segment lies in the
        interior of a ball of that orbit with radius reduced by \(\delta\).
        Each such set \(U\) is open. These sets cover
        \(\mathbb T\times S^2\), since every segment of length \(L\) enters
        an original ball's interior and retains an interior intersection
        after a sufficiently small positive radius reduction.

        A finite subcover provides finitely many
        positive reductions. Choose \(\eta\) smaller than all of them and
        smaller than every original radius. Every segment of length \(L\) then
        enters a ball with radius reduced by \(\eta\).
    \end{proof}

    \begin{corollary}[Stability under independent perturbations]
        \label{cor:stability}
        Index the balls in \(\mathcal B\) by \(\alpha\), with centers
        \(\vect{c}_\alpha\) and radii \(r_\alpha\). There is an
        \(\varepsilon>0\) such that any family of positive-radius balls satisfying
        \begin{equation}
            \sup_\alpha\bigl(
            \norm{\vect{c}'_\alpha-\vect{c}_\alpha}
            +\abs{r'_\alpha-r_\alpha}\bigr)<\varepsilon
            \label{eq:stability-perturbation}
        \end{equation}
        is strictly separated and has finite horizon.
    \end{corollary}

    \begin{proof}

        By periodicity and strict separation, choose \(g>0\) such that
        \[
            \norm{\vect{c}_\alpha-\vect{c}_\beta}
            -(r_\alpha+r_\beta)\ge g
            \qquad\text{for every }\alpha\ne\beta.
        \]
        Thus \(g\) is a uniform lower bound on the distances between
        distinct ball boundaries.

        By Proposition~\ref{prop:horizon-upper} and
        Lemma~\ref{lem:uniform-segment-margin}, there are \(L,\eta>0\) such
        that every segment of length \(L\) enters a witness ball of radius
        \(r_\alpha-\eta\) centered at \(\vect{c}_\alpha\).

        Choose \(0<\varepsilon<\min\{g/2,\eta\}\).
        By~\eqref{eq:stability-perturbation}, the gap between any two
        perturbed balls is greater than \(g-2\varepsilon>0\). The same bound
        and \(\varepsilon<\eta\) ensure that each perturbed ball indexed by
        \(\alpha\) contains the witness ball centered at \(\vect{c}_\alpha\)
        in its interior. Hence every segment of length \(L\) enters a perturbed
        ball.
    \end{proof}

    \section{Concluding remarks}

    The construction uses 38 scatterers per period cell, and
    the choices \(\BlockingMargin=1/250\) and \(N=360\) simplify the estimates
    leading to \(\Horizon<453601\). The minimum number of scatterers needed
    per period cell for a strictly separated three-dimensional configuration
    with finite horizon remains to be determined, as does the horizon of the
    present example.

    Corollary~\ref{cor:stability} shows that strict separation
    and finite horizon persist under uniformly small independent perturbations
    of the centers and radii. The construction uses four radii, leaving open
    Dettmann's question of finding an explicit periodic arrangement of
    equal-sized, nonoverlapping spheres with finite horizon in dimension
    three~\cite[Open Problem~2]{Dettmann14}.

    \appendix

    \section{Lattice facts and planar covering formulas}
    \label{app:lattice-facts}

    This appendix supplies the lattice facts used in
    Section~\ref{sec:rational-blocking}. The main text uses their conclusions
    to estimate projected covering radii; the derivations are collected here
    in their order of dependence.

    \subsection{Integer bases and orthogonal projection}
    \label{app:lattice-bases}

    Basis extension and Gauss reduction use the following integer changes
    of basis.
    Let \(\matr{A}\) be a basis matrix for a lattice
    \(\latt{L}=\matr{A}\Z^n\).
    If an invertible \(n\times n\) matrix \(V\) and its inverse
    \(V^{-1}\) have integer entries, then
    \begin{equation}
        V\Z^n=\Z^n,
        \qquad (\matr{A}V)\Z^n
        =\matr{A}\Z^n=\latt{L}.
        \label{eq:integer-basis-change}
    \end{equation}
    Thus the columns of \(\matr{A}V\) form another basis of \(\latt{L}\).
    Swapping basis vectors, changing signs, and subtracting integer multiples
    of one basis vector from another correspond to matrices \(V\) of this kind.

    A primitive vector \(\vect{w}\in\latt{L}\) extends to a basis of
    \(\latt{L}\). To see this, choose a lattice basis matrix
    \(\matr{A}\) and write
    \(\vect{w}=\matr{A}\vect{a}\), where \(\vect{a}\in\Z^n\) has
    coordinates of greatest common divisor one. The Euclidean algorithm
    gives a matrix \(U\) with integer entries and an integer inverse such
    that \(U\vect{a}=\vect{e}_1\). By~\eqref{eq:integer-basis-change},
    the columns of \(\matr{A}U^{-1}\) form a basis of \(\latt{L}\), and the
    first column is
    \(\matr{A}U^{-1}\vect{e}_1
    =\matr{A}\vect{a}=\vect{w}\).
    The condition \(\latt{L}\cap\R\vect{w}=\Z\vect{w}\) also shows
    that primitivity is independent of the chosen lattice basis.

    In~\eqref{eq:orthogonal-projection}, the component
    along \(\hat{\vect{w}}\) is orthogonal to the projected component.
    Thus for \(\vect{u},\vect{v}\in\R^3\),
    \begin{equation}
        \orthproj{w}(\vect{u})\dotprod\orthproj{w}(\vect{v})
        =\vect{u}\dotprod\vect{v}
        -(\vect{u}\dotprod\hat{\vect{w}})(\vect{v}\dotprod\hat{\vect{w}}).
        \label{eq:projected-inner-product}
    \end{equation}
    Taking \(\vect{v}=\vect{u}\) gives
    \(\norm{\orthproj{w}(\vect{u})}\le\norm{\vect{u}}\).

    \subsection{Projected covolumes and short vectors}

    \begin{lemma}[Covolume after projection]
        \label{lem:projected-covolume}
        Let \(\latt{L}\subset\R^3\) be a full-rank lattice and let
        \(\vect{w}\in\latt{L}\) be primitive. Then
        \[
            \covol\bigl(\orthproj{w}(\latt{L})\bigr)
            =\frac{\covol(\latt{L})}{\norm{\vect{w}}}.
        \]
        In particular,
        \(\covol(\projlat{w})=2/\norm{\vect{w}}\) for \(\latt{L}=\FCC\).
    \end{lemma}

    \begin{proof}
        Complete the primitive vector \(\vect{w}\) to a lattice basis
        \(\vect{w},\vect{b}_1,\vect{b}_2\) of \(\latt{L}\). The vectors
        \(\orthproj{w}(\vect{b}_1)\) and
        \(\orthproj{w}(\vect{b}_2)\) form a basis of the projected lattice.
        If their fundamental parallelogram has area \(\Area\), the base-times-height
        formula gives
        \[
            \covol(\latt{L})=\norm{\vect{w}}\Area.
        \]
        Since \(\Area=\covol(\orthproj{w}(\latt{L}))\), division by
        \(\norm{\vect{w}}\) proves the formula. Primitivity ensures that the
        height is one lattice period.
    \end{proof}

    We use planar basis reduction both to find a short projected vector and
    to choose bases for the exact covering-radius calculations.

    \begin{lemma}[Gauss reduction]
        \label{lem:gauss-reduction}
        Any basis of a planar lattice \(\Gamma\) can be transformed in
        finitely many steps, by swaps, subtraction of integer multiples of
        one basis vector from the other, and sign changes, into a basis
        \(\planarvec{\xi}_1,\planarvec{\xi}_2\) satisfying
        \begin{equation}
            \norm{\planarvec{\xi}_1}\le\norm{\planarvec{\xi}_2},
            \qquad
            -\frac12\norm{\planarvec{\xi}_1}^2
            \le\planarvec{\xi}_1\dotprod\planarvec{\xi}_2\le0.
            \label{eq:gauss-reduced-basis}
        \end{equation}
    \end{lemma}

    \begin{proof}
        Each basis change preserves \(\Gamma\) by
        \eqref{eq:integer-basis-change}. First order the basis vectors so that
        \(\norm{\planarvec{\xi}_1}\le\norm{\planarvec{\xi}_2}\).
        Choose \(k\in\Z\) as a nearest integer to
        \[
            \frac{\planarvec{\xi}_1\dotprod\planarvec{\xi}_2}
            {\norm{\planarvec{\xi}_1}^2},
        \]
        and replace \(\planarvec{\xi}_2\) by
        \(\planarvec{\xi}_2-k\planarvec{\xi}_1\).
        Subtracting \(k\planarvec{\xi}_1\) shifts the projection of
        \(\planarvec{\xi}_2\) onto \(\R\planarvec{\xi}_1\) by integer
        steps of length \(\norm{\planarvec{\xi}_1}\). The choice of
        \(k\) places the projected vector between \(-\planarvec{\xi}_1/2\) and
        \(\planarvec{\xi}_1/2\), giving
        \[
            \abs{\planarvec{\xi}_1\dotprod\planarvec{\xi}_2}
            \le\frac12\norm{\planarvec{\xi}_1}^2.
        \]

        If the new \(\planarvec{\xi}_2\) is shorter than
        \(\planarvec{\xi}_1\), swap the basis vectors and repeat the
        subtraction. Every such swap strictly decreases
        \(\norm{\planarvec{\xi}_1}\). A lattice has only finitely many
        vectors in a bounded disk, so only finitely many swaps can occur.
        When the new \(\planarvec{\xi}_2\) is at least as long as
        \(\planarvec{\xi}_1\), stop. Changing the sign of
        \(\planarvec{\xi}_2\), if necessary, makes the inner product
        nonpositive and gives~\eqref{eq:gauss-reduced-basis}.
    \end{proof}

    \begin{lemma}[A short vector in a planar lattice]
        \label{lem:planar-short-vector}
        Every two-dimensional lattice of covolume \(\Area\) contains a nonzero
        vector \(\planarvec{\xi}\) such that
        \[
            \norm{\planarvec{\xi}}^2\le\frac{2\Area}{\sqrt3}.
        \]
    \end{lemma}

    \begin{proof}
        Denote the planar lattice by \(\Gamma\), and choose a basis
        \(\planarvec{\xi}_1,\planarvec{\xi}_2\) satisfying
        \eqref{eq:gauss-reduced-basis} by
        Lemma~\ref{lem:gauss-reduction}. The angle \(\theta\) between
        the basis vectors satisfies \(-1/2\le\cos\theta\le0\), so
        \(\pi/2\le\theta\le2\pi/3\). Therefore
        \[
            \Area
            =\norm{\planarvec{\xi}_1}\norm{\planarvec{\xi}_2}\sin\theta
            \ge\frac{\sqrt3}{2}\norm{\planarvec{\xi}_1}^2,
        \]
        giving the claimed bound with \(\planarvec{\xi}=\planarvec{\xi}_1\).
    \end{proof}

    Lemmas~\ref{lem:projected-covolume}
    and~\ref{lem:planar-short-vector} give a nonzero vector
    \(\planarvec{\xi}\in\orthproj{w}(\latt{L})\) with
    \[
        \norm{\planarvec{\xi}}^2
        \leq\frac{2\covol(\latt{L})}
        {\sqrt3\,\norm{\vect{w}}}.
    \]
    For \(\latt{L}=\FCC\), this reads
    \begin{equation}
        \norm{\planarvec{\xi}}^2
        \le\frac4{\sqrt3\,\norm{\vect{w}}}.
        \label{eq:short-projection}
    \end{equation}

    \subsection{Dual lattices and plane spacing}

    We next derive the plane-spacing bound from the dual lattice defined in
    Section~\ref{subsec:projected-estimate}.
    If \(\latt{L}=\matr{A}\Z^n\), then
    \[
        \vect{y}\in\latt{L}^*
        \iff \matr{A}^T\vect{y}\in\Z^n.
    \]
    Therefore
    \begin{equation}
        \latt{L}^*=\matr{A}^{-T}\Z^n,
        \label{eq:dual-basis-formula}
    \end{equation}
    and the columns of \(\matr{A}^{-T}\)
    form the dual basis corresponding to the columns of \(\matr{A}\).

    Let \(\vect{y}\in\latt{L}^*\) be primitive. By the definition of
    the dual lattice, \(\vect{y}\dotprod\latt{L}\subseteq\Z\).
    Since the coordinates of \(\vect{y}\) in a dual basis have
    greatest common divisor one, B\'ezout's identity gives
    \(\vect{y}\dotprod\latt{L}\supseteq\Z\).
    Thus \(\vect{y}\dotprod\latt{L}=\Z\), so every hyperplane
    \begin{equation}
        H_k=\{\vect{x}\in\R^n:\vect{y}\dotprod\vect{x}=k\},
        \qquad k\in\Z,
        \label{eq:dual-lattice-layers}
    \end{equation}
    contains a lattice point. Consecutive hyperplanes \(H_k\) and
    \(H_{k+1}\) are spaced by \(1/\norm{\vect{y}}\).

    To apply~\eqref{eq:dual-lattice-layers} to a prescribed plane
    \(V\subset\R^3\) spanned by lattice vectors, we need a primitive
    dual vector perpendicular to \(V\). Write
    \(\latt{L}=\matr{A}\Z^3\) and
    \(V=\Span(\matr{A}\vect{u},\matr{A}\vect{v})\), where
    \(\vect{u},\vect{v}\in\Z^3\) are linearly independent. Set
    \begin{equation}
        \vect{c}=\vect{u}\times\vect{v},\qquad
        d=\gcd(c_1,c_2,c_3),\qquad
        \vect{y}=\matr{A}^{-T}\frac{\vect{c}}{d}.
        \label{eq:primitive-dual-normal}
    \end{equation}
    Since \(\latt{L}^*=\matr{A}^{-T}\Z^3\) and the coordinates of
    \(\vect{c}/d\) have greatest common divisor one,
    \(\vect{y}\) is primitive in \(\latt{L}^*\). Moreover,
    \[
        \vect{y}\dotprod \matr{A}\vect{u}
        =\frac{\vect{c}\dotprod\vect{u}}{d}=0,
        \qquad
        \vect{y}\dotprod \matr{A}\vect{v}
        =\frac{\vect{c}\dotprod\vect{v}}{d}=0,
    \]
    so \(\vect{y}\in V^\perp\).

    \begin{lemma}[Lattice plane-spacing estimate]
        \label{lem:lattice-plane-spacing}
        Let \(\latt{L}\subset\R^3\) be a full-rank lattice, and let \(V\) be
        spanned by two linearly independent vectors of \(\latt{L}\). Every affine
        plane parallel to \(V\) lies at distance at most
        \[
            \frac{1}{2\lambda_1(\latt{L}^*)}
        \]
        from a parallel plane containing lattice points. Hence, if
        \(R>1/(2\lambda_1(\latt{L}^*))\), every such affine plane meets the
        interior of an \(R\)-ball centered in \(\latt{L}\).
    \end{lemma}

    \begin{proof}
        Choose a primitive dual normal \(\vect{y}\in\latt{L}^*\cap V^\perp\)
        by~\eqref{eq:primitive-dual-normal}. The planes \(H_k\)
        in~\eqref{eq:dual-lattice-layers} are parallel to \(V\), contain
        lattice points, and have spacing \(1/\norm{\vect{y}}\).
        Given an affine plane \(P\parallel V\), choose a nearest plane
        \(H_k\). Then
        \[
            \dist(H_k,P)\le\frac{1}{2\norm{\vect{y}}}
            \le\frac{1}{2\lambda_1(\latt{L}^*)}.
        \]
        Choose \(\vect{p}\in\latt{L}\cap H_k\). If
        \(R>1/(2\lambda_1(\latt{L}^*))\), then
        \(\dist(\vect{p},P)=\dist(H_k,P)<R\), so \(P\) meets the interior
        of the radius-\(R\) ball centered at \(\vect{p}\).
    \end{proof}

    \subsection{The FCC dual lattice}

    \begin{lemma}[FCC dual lattice and minimum]
        \label{lem:fcc-dual}
        The dual of the FCC lattice is
        \[
            \dualFCC=\Z^3\cup
            \bigl((1/2,1/2,1/2)+\Z^3\bigr),
        \]
        and its shortest nonzero vectors have length \(\sqrt3/2\).
    \end{lemma}

    \begin{proof}
        The vectors
        \[
            \vect{a}_1=(1,1,0),\qquad \vect{a}_2=(1,0,1),\qquad
            \vect{a}_3=(0,1,1)
        \]
        form a basis of \(\FCC\). Set
        \(\matr{A}=(\vect{a}_1\mid\vect{a}_2\mid\vect{a}_3)\).
        By~\eqref{eq:dual-basis-formula},
        \[
            \dualFCC=\matr{A}^{-T}\Z^3,\qquad
            \matr{A}^{-T}=\frac12
            \begin{pmatrix}
                1&1&-1\\
                1&-1&1\\
                -1&1&1
            \end{pmatrix}.
        \]
        Each column of \(\matr{A}^{-T}\) has all three coordinates congruent
        to \(1/2\) modulo \(\Z\). Every integer combination therefore
        has either three integer coordinates or three half-integer
        coordinates, proving
        \[
            \matr{A}^{-T}\Z^3
            \subseteq\Z^3\cup\bigl((1/2,1/2,1/2)+\Z^3\bigr).
        \]
        For the reverse inclusion, the pairwise sums of the columns
        of \(\matr{A}^{-T}\) are the standard basis vectors
        \(\vect{e}_1,\vect{e}_2,\vect{e}_3\), so
        \(\Z^3\subseteq \matr{A}^{-T}\Z^3\).
        The sum of all three columns is \((1/2,1/2,1/2)\), so
        \((1/2,1/2,1/2)+\Z^3\subseteq \matr{A}^{-T}\Z^3\) as well.

        A nonzero vector in the integer coset has length at least \(1\),
        whereas a vector in the half-integer coset has length at least
        \(\sqrt3/2\). The latter bound is attained by \((1/2,1/2,1/2)\),
        so \(\lambda_1(\dualFCC)=\sqrt3/2\).
    \end{proof}

    \subsection{Exact planar covering radii}
    \label{app:planar-covering}

    The planar form of Selling--Voronoi reduction gives the covering formula
    used for the six short FCC types. The standard terminology and general
    theory are given in \cite[Theorem~3]{CS92}; we include the planar proof.

    A \emph{superbase} of a planar lattice \(\Gamma\) is a triple
    \(\planarvec{\xi}_0,\planarvec{\xi}_1,\planarvec{\xi}_2\)
    whose sum is zero and any two of which form a basis of \(\Gamma\).
    A superbase is \emph{obtuse} if
    \(\planarvec{\xi}_i\dotprod\planarvec{\xi}_j\le0\)
    whenever \(i\neq j\).

    We first show how obtuseness excludes lattice points from the interiors
    of the circumdisks of the two triangles in a fundamental parallelogram.

    \begin{lemma}[Empty circumdisks from an obtuse superbase]
        \label{lem:obtuse-empty-circumdisks}
        Let \(\planarvec{\xi}_0,\planarvec{\xi}_1,\planarvec{\xi}_2\)
        be an obtuse superbase of a planar lattice \(\Gamma\).
        Let \(T_1\) have vertices
        \(0,\planarvec{\xi}_1,\planarvec{\xi}_1+\planarvec{\xi}_2\),
        and let \(T_2\) have vertices
        \(0,\planarvec{\xi}_2,\planarvec{\xi}_1+\planarvec{\xi}_2\).
        Neither circumdisk contains a point of \(\Gamma\) in its interior.
    \end{lemma}

\begin{figure}[htbp]
        \centering
        \includegraphics[width=\linewidth]{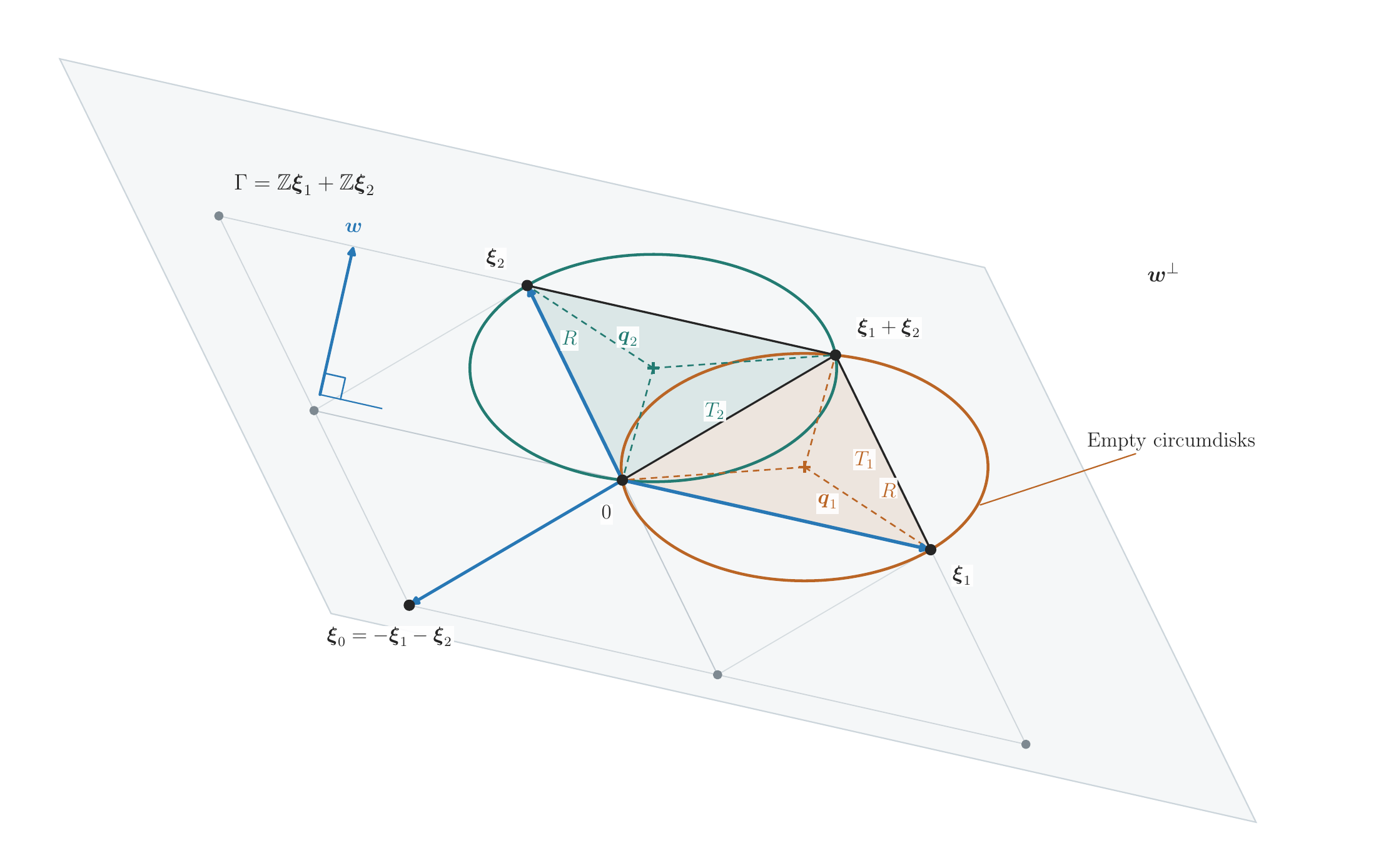}
        \caption{The two triangles of the basis parallelogram and their empty
            circumdisks, shown in the projection plane \(\orthcomp{w}\) of the
            FCC application. Their common circumradius \(R\) is the covering
            radius of \(\Gamma\).}
        \label{fig:obtuse-superbase}
    \end{figure}

    Figure~\ref{fig:obtuse-superbase} shows the two triangles and their circumdisks.

    \begin{proof}
        Let \(\planarvec{q}_1\) be the circumcenter of \(T_1\).
        Since \(0\) is a vertex of \(T_1\), the circumradius of \(T_1\)
        is \(\norm{\planarvec{q}_1}\). To exclude lattice points from the
        interior of the circumdisk of \(T_1\), we will prove
        \[
            \norm{\planarvec{z}-\planarvec{q}_1}^2
            \ge\norm{\planarvec{q}_1}^2
            \qquad\text{for every }\planarvec{z}\in\Gamma.
        \]

        Write \(g_{ij}=\planarvec{\xi}_i\dotprod\planarvec{\xi}_j\)
        for \(i,j\in\{1,2\}\), and set
        \[
            \begin{aligned}
                p_1&=-\planarvec{\xi}_0\dotprod\planarvec{\xi}_1
                    =g_{11}+g_{12},\\
                p_2&=-\planarvec{\xi}_0\dotprod\planarvec{\xi}_2
                    =g_{22}+g_{12},\\
                p_3&=-\planarvec{\xi}_1\dotprod\planarvec{\xi}_2
                    =-g_{12}.
            \end{aligned}
        \]
        The equalities use
        \(\planarvec{\xi}_0=-\planarvec{\xi}_1-\planarvec{\xi}_2\),
        and obtuseness gives \(p_1,p_2,p_3\ge0\).

        Equidistance from the three vertices gives
        \[
            \norm{\planarvec{q}_1}^2
            =\norm{\planarvec{q}_1-\planarvec{\xi}_1}^2
            =\norm{\planarvec{q}_1-\planarvec{\xi}_1-\planarvec{\xi}_2}^2.
        \]
        Expanding the last two squared norms yields
        \begin{align*}
            0&=g_{11}-2\planarvec{q}_1\dotprod\planarvec{\xi}_1,\\
            0&=g_{11}+2g_{12}+g_{22}
                -2\planarvec{q}_1\dotprod\planarvec{\xi}_1
                -2\planarvec{q}_1\dotprod\planarvec{\xi}_2.
        \end{align*}
        Consequently,
        \[
            2\planarvec{q}_1\dotprod\planarvec{\xi}_1=g_{11},
            \qquad
            2\planarvec{q}_1\dotprod\planarvec{\xi}_2=g_{22}+2g_{12}.
        \]

        Every lattice point has the form
        \(\planarvec{z}=m\planarvec{\xi}_1+n\planarvec{\xi}_2\)
        with \(m,n\in\Z\). Using the circumcenter identities above and
        \(g_{11}=p_1+p_3\), \(g_{22}=p_2+p_3\), \(g_{12}=-p_3\),
        we obtain
        \begin{align*}
            \norm{\planarvec{z}-\planarvec{q}_1}^2
                -\norm{\planarvec{q}_1}^2
            &=g_{11}m(m-1)+g_{22}n(n-1)\\
            &\quad+2g_{12}n(m-1)\\
            &=p_1m(m-1)+p_2n(n-1)\\
            &\quad+p_3(m-n)(m-n-1)\ge0.
        \end{align*}
        The last inequality follows because obtuseness gives
        \(p_1,p_2,p_3\ge0\), while \(k(k-1)\ge0\) for every integer \(k\).
        Thus no lattice point lies in the interior
        of the circumdisk of \(T_1\). Interchanging
        \(\planarvec{\xi}_1\) and \(\planarvec{\xi}_2\) proves the
        assertion for \(T_2\).
    \end{proof}

    Lemma~\ref{lem:obtuse-empty-circumdisks} identifies disks with no lattice
    points in their interiors. We now use the triangles \(T_1,T_2\) to show
    that their common circumradius is the covering radius and compute its value.

    \begin{lemma}[Planar covering formula for an obtuse superbase]
        \label{lem:planar-covering}
        Let \(\planarvec{\xi}_1,\planarvec{\xi}_2\) be linearly independent
        vectors in a Euclidean plane, set
        \(\Gamma=\Z\planarvec{\xi}_1+\Z\planarvec{\xi}_2\), and put
        \[
            g_{11}=\norm{\planarvec{\xi}_1}^2,\qquad
            g_{12}=\planarvec{\xi}_1\dotprod\planarvec{\xi}_2,\qquad
            g_{22}=\norm{\planarvec{\xi}_2}^2.
        \]
        Linear independence gives \(g_{11}g_{22}-g_{12}^2>0\).
        If \(g_{12}\leq0\), \(g_{11}+g_{12}\geq0\), and
        \(g_{22}+g_{12}\geq0\), then, with
        \(\planarvec{\xi}_0=-\planarvec{\xi}_1-\planarvec{\xi}_2\), the vectors
        \(\planarvec{\xi}_0,\planarvec{\xi}_1,\planarvec{\xi}_2\) form an obtuse
        superbase, and
        \begin{equation}
            \mu(\Gamma)^2=\frac{g_{11}g_{22}(g_{11}+g_{22}+2g_{12})}
            {4(g_{11}g_{22}-g_{12}^2)}.
            \label{eq:exact-short-covering-radius}
        \end{equation}
    \end{lemma}

    \begin{proof}
        Since \(\planarvec{\xi}_0=-\planarvec{\xi}_1-\planarvec{\xi}_2\),
        any two of \(\planarvec{\xi}_0,\planarvec{\xi}_1,\planarvec{\xi}_2\)
        form a basis of \(\Gamma\). The identities
        \[
            \planarvec{\xi}_0\dotprod\planarvec{\xi}_1=-(g_{11}+g_{12}),
            \qquad
            \planarvec{\xi}_0\dotprod\planarvec{\xi}_2=-(g_{22}+g_{12})
        \]
        and \(\planarvec{\xi}_1\dotprod\planarvec{\xi}_2=g_{12}\)
        show that the superbase is obtuse under the stated hypotheses.
        Let \(T_1,T_2\) be the triangles in
        Lemma~\ref{lem:obtuse-empty-circumdisks}, and let
        \(\planarvec{q}_1\) be the circumcenter of \(T_1\).
        The circumdisks of \(T_1,T_2\) therefore contain no lattice points
        in their interiors.

        The triangles \(T_1,T_2\) have the same side lengths and hence
        the same circumradius \(R\). Their lattice translates tile the plane.
        To show that \(R\) is the covering radius, consider any triangle
        \(T\) in this tiling, with vertices
        \(\planarvec{v}_0,\planarvec{v}_1,\planarvec{v}_2\)
        and circumcenter \(\planarvec{q}_T\).
        Every point \(\planarvec{x}\in T\) has a representation
        \[
            \planarvec{x}=\sum_{i=0}^2\lambda_i\planarvec{v}_i,
            \qquad \lambda_i\ge0,\qquad \sum_{i=0}^2\lambda_i=1.
        \]
        Since \(\norm{\planarvec{v}_i-\planarvec{q}_T}=R\),
        \[
            \sum_{i=0}^2\lambda_i
            \norm{\planarvec{x}-\planarvec{v}_i}^2
            =R^2-\norm{\planarvec{x}-\planarvec{q}_T}^2
            \le R^2.
        \]
        At least one vertex of \(T\) is therefore within \(R\) of
        \(\planarvec{x}\), proving \(\mu(\Gamma)\le R\).
        The empty circumdisk of \(T_1\) gives
        \(\dist(\planarvec{q}_1,\Gamma)=R\), so
        \(\mu(\Gamma)=R\).

        Finally, the squared side lengths of \(T_1\) are
        \(g_{11},g_{22},g_{11}+g_{22}+2g_{12}\), and its squared area is
        \((g_{11}g_{22}-g_{12}^2)/4\).
        Substituting these values into the circumradius formula
        \(R=abc/(4\Delta)\), where \(a,b,c\) are the side lengths
        and \(\Delta\) is the area, gives
        \eqref{eq:exact-short-covering-radius}.
    \end{proof}

    For each representative \(\vect{w}\) in
    Table~\ref{tab:short-fcc-radii}, complete \(\vect{w}\) to an FCC basis
    \(\vect{w},\vect{b}_1,\vect{b}_2\) using the basis-extension property
    in Appendix~\ref{app:lattice-bases}. Put
    \(\planarvec{\xi}_i=\orthproj{w}(\vect{b}_i)\).
    Apply Lemma~\ref{lem:gauss-reduction} to the projected basis
    \(\planarvec{\xi}_1,\planarvec{\xi}_2\), performing the same swaps,
    integer subtractions, and sign changes on \(\vect{b}_1,\vect{b}_2\).
    The operations on \(\vect{b}_1,\vect{b}_2\) preserve the FCC basis by
    \eqref{eq:integer-basis-change} and give projected vectors
    satisfying~\eqref{eq:gauss-reduced-basis}. In terms of the Gram triple
    \[
        (g_{11},g_{12},g_{22})
        =\bigl(\norm{\planarvec{\xi}_1}^2,
        \planarvec{\xi}_1\dotprod\planarvec{\xi}_2,
        \norm{\planarvec{\xi}_2}^2\bigr),
    \]
    the inequalities are \(g_{11}\le g_{22}\) and
    \(-g_{11}/2\le g_{12}\le0\), which imply the hypotheses of
    Lemma~\ref{lem:planar-covering}.
    Table~\ref{tab:short-fcc-radii} lists suitable
    basis choices and the Gram triples obtained by applying
    \eqref{eq:projected-inner-product} to the pairs \((\vect{b}_i,\vect{b}_j)\).

    \section{Exact line-space certificate}
    \label{app:certificate}

    We derive the cell-containment test used in Lemma~\ref{lem:certificate}.
    Let \(\planarvec{p}\) be the center of a planar cell,
    and suppose every point in the cell is within distance \(\rho\) of
    \(\planarvec{p}\). By the triangle inequality, the open disk of radius \(R\)
    centered at \(\planarvec{q}\) contains the cell whenever
    \begin{equation}
        \norm{\planarvec{p}-\planarvec{q}}+\rho<R.
        \label{eq:cell-containment-test}
    \end{equation}
    We use the dyadic cell center~\eqref{eq:dyadic-cell-center}
    and the cell distance bound~\eqref{eq:cell-radius} to express the
    witness-ball containment test in integer arithmetic.

    Fix a primitive vector \(\vect{a}\in\Z^3\). The basis-extension property
    in Appendix~\ref{app:lattice-bases} gives
    \(\vect{b}_1,\vect{b}_2\in\Z^3\) such that
    \[
        \det[\vect{a}\ \vect{b}_1\ \vect{b}_2]=\pm1.
    \]

    As in Section~\ref{subsec:primitive-periods}, project
    the basis \((2\vect{a},2\vect{b}_1,2\vect{b}_2)\) of \(2\Z^3\) onto
    \(\orthcomp{a}\). The vectors \(2\orthproj{a}(\vect{b}_1)\) and
    \(2\orthproj{a}(\vect{b}_2)\) form a basis of \(\orthproj{a}(2\Z^3)\).
    Since \(\vect{a}\) is parallel to
    \(\vect{w}(\vect{a})\), the resulting fundamental parallelogram gives
    coordinates \((s,t)\) on the line space \(\mathcal S_{\vect{w}(\vect{a})}\).
    Each line class has a unique representative of the form
    \[
        \ell_{s,t}=2(s\vect{b}_1+t\vect{b}_2)+\R\vect{a},
        \qquad (s,t)\in[0,1)^2.
    \]
    Appendix~\ref{alg:line-certificate} specifies
    \(\vect{b}_1,\vect{b}_2\) by reducing the projected basis.

    At depth \(d\), set \(N=2^d\) and divide the parameter
    square \([0,1)^2\) into \(N^2\) congruent dyadic squares. The midpoint of
    the square with indices \(0\leq i,j<N\) represents the line through
    the base point
    \begin{equation}
        \vect{x}_{ij}=\frac{(2i+1)\vect{b}_1+(2j+1)\vect{b}_2}{N}.
        \label{eq:dyadic-cell-center}
    \end{equation}

    For a vector \(\vect{x}=(x_1,x_2,x_3)\), write
    \[
        \onenorm{\vect{x}}=\abs{x_1}+\abs{x_2}+\abs{x_3},
        \qquad
        L=\onenorm{\vect{b}_1}+\onenorm{\vect{b}_2}.
    \]
    If a parameter belongs to this dyadic square, the corresponding base point
    is at distance at most
    \begin{equation}
        \rho_N=\frac{L}{N}
        \label{eq:cell-radius}
    \end{equation}
    from \(\vect{x}_{ij}\). The two parameter displacements are at most
    \(1/(2N)\), and the factor \(2\) in \(\ell_{s,t}\)
    gives~\eqref{eq:cell-radius} by the triangle inequality. The one-norm keeps
    this bound rational.
    By~\eqref{eq:projected-inner-product}, every point of
    the projected cell is also within \(\rho_N\) of
    \(\orthproj{a}(\vect{x}_{ij})\).

    The bound \(\rho_N\) accounts for the size of the
    projected cell. To apply the containment test~\eqref{eq:cell-containment-test},
    we now compute the squared distance from \(\orthproj{a}(\vect{x}_{ij})\)
    to the center of a candidate witness disk.

    Choose a scatterer of radius \(r=U/100\), \(U\in\Z\),
    with center \(\vect{c}+2\vect{k}\), where \(\vect{c}\) is a listed center
    and \(\vect{k}\in\Z^3\). Let \(\ell_{ij}=\vect{x}_{ij}+\R\vect{a}\) be the
    midpoint line and set \(s_{ij}=\dist(\vect{c}+2\vect{k},\ell_{ij})\).
    Projecting the displacement \(\vect{x}_{ij}-(\vect{c}+2\vect{k})\) onto
    \(\orthcomp{a}\) gives a vector of length \(s_{ij}\).

\begin{figure}[htbp]
        \centering
        \includegraphics[width=\linewidth]{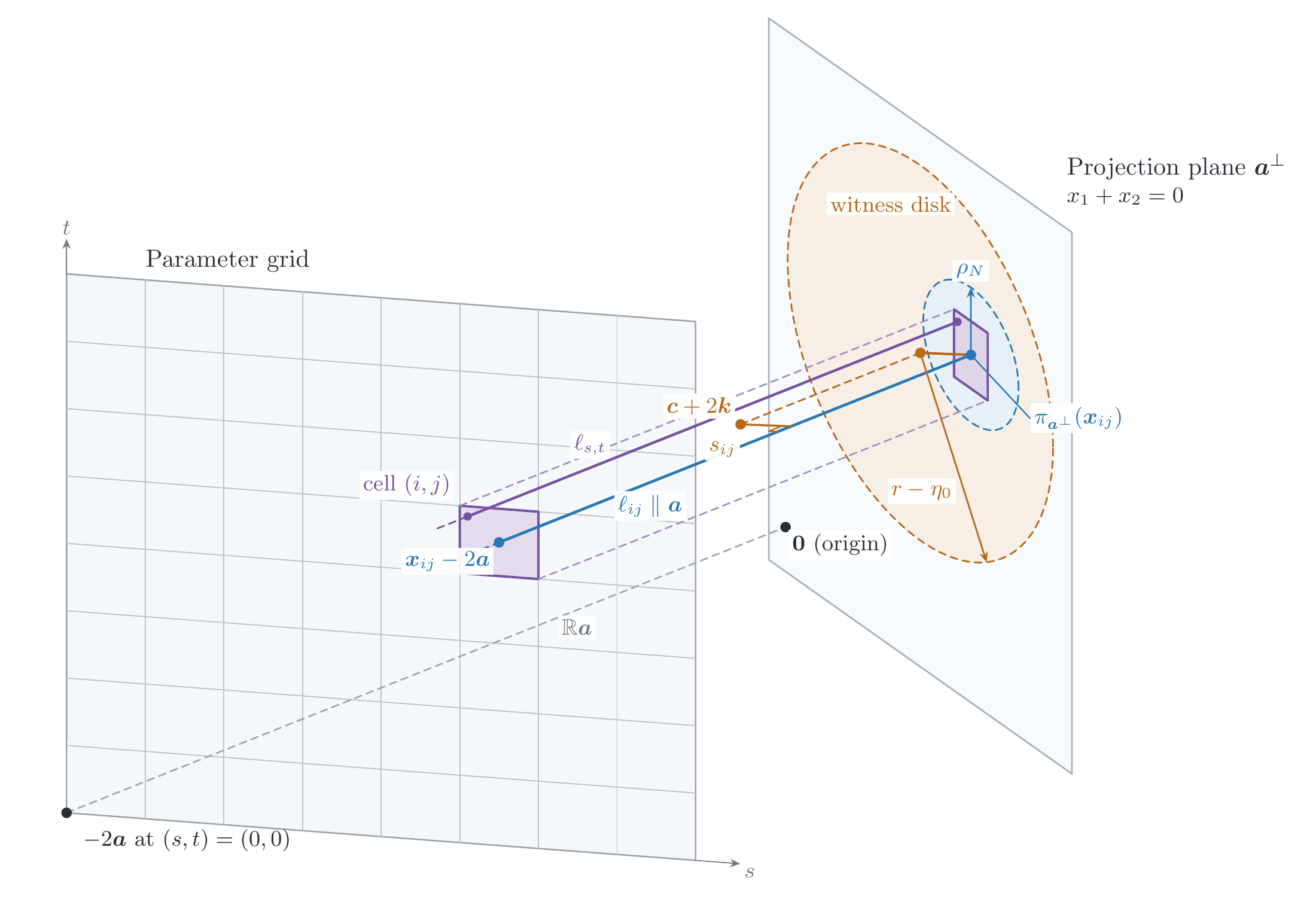}
        \caption{The cell test for \(\vect{a}=(1,1,0)\). The parameter grid
            is drawn at base points \(2(s\vect{b}_1+t\vect{b}_2)-2\vect{a}\)
            on the same lines \(\ell_{s,t}\).
            The projected cell lies inside the witness disk when
            \(s_{ij}+\rho_N<r-\eta_0\).}
        \label{fig:line-space-certificate}
    \end{figure}

    Figure~\ref{fig:line-space-certificate} shows the cell, its corresponding
    lines, and their projection.

    To compute \(s_{ij}^2\) using integers, first clear
    the coordinate denominators in the displacement. The factor \(N\)
    clears the denominators of \(\vect{x}_{ij}\) by~\eqref{eq:dyadic-cell-center},
    and \(10\) clears the denominators of \(\vect{c}+2\vect{k}\) because the
    listed center coordinates lie in \(\tfrac1{10}\Z\).
    With \(\vect{C}=10\vect{c}\) and \(\vect{Z}=\vect{C}+20\vect{k}\),
    the scaled displacement is
    \[
        \vect{D}=10N\bigl(\vect{x}_{ij}-(\vect{c}+2\vect{k})\bigr)
        =10N\vect{x}_{ij}-N\vect{Z}\in\Z^3.
    \]
    Orthogonal projection is linear, so the factor
    \(10N\) also scales the projected displacement. The vector
    \(\orthproj{a}(\vect{D})\) therefore has length \(10N s_{ij}\).
    Setting \(\vect{u}=\vect{v}=\vect{D}\)
    in~\eqref{eq:projected-inner-product} gives
    \begin{equation}
        s_{ij}^2
        =\frac{\norm{\orthproj{a}(\vect{D})}^2}{(10N)^2}
        =\frac{\norm{\vect{a}}^2\norm{\vect{D}}^2
        -(\vect{D}\dotprod\vect{a})^2}
        {100N^2\norm{\vect{a}}^2}.
        \label{eq:projected-distance}
    \end{equation}

    \begin{claim}[Dyadic cell acceptance]
        The entire projected cell with center
        \(\orthproj{a}(\vect{x}_{ij})\) lies in the open projected
        disk of the witness ball of radius \((5U-2)/500\) centered at
        \(\vect{c}+2\vect{k}\) whenever
        \((5U-2)N-500L>0\) and
        \begin{equation}
            2500\bigl(\norm{\vect{a}}^2\norm{\vect{D}}^2
            -(\vect{D}\dotprod\vect{a})^2\bigr)
            <\norm{\vect{a}}^2\bigl((5U-2)N-500L\bigr)^2.
            \label{eq:reduced-radius-integer-certificate}
        \end{equation}
    \end{claim}

    \begin{proof}
        For the projected cell,
        \eqref{eq:cell-containment-test} asks for \(s_{ij}+\rho_N<r-\BlockingMargin\).
        Since \(r-\BlockingMargin=(5U-2)/500\) and \(\rho_N=L/N\), the first
        inequality in the claim ensures \(r-\BlockingMargin-\rho_N>0\).
        We can therefore square the comparison with \(s_{ij}\). Substituting
        \eqref{eq:projected-distance} and clearing denominators gives exactly
        the integer inequality~\eqref{eq:reduced-radius-integer-certificate}.
    \end{proof}

    For each of the nine exceptional directions in
    Lemma~\ref{lem:directions}, the algorithm in Appendix~\ref{alg:line-certificate}
    uses the dyadic cell test to certify that every line represented by the
    parameter square enters a witness ball.

    \section{Exact verification algorithms}
    \label{app:reproducibility}

    The algorithms below use the center and radius data in
    Section~\ref{sec:explicit}. In particular, the line-space procedure generates
    and tests its dyadic cells rather than reading a precomputed certificate.

    The computations produce the following exact values:
    \begin{center}
        \begin{tabular}{lc}
            \toprule
            quantity                                             & exact value      \\
            \midrule
            scatterers on \(\PeriodTorus\)                       & 38               \\
            short directions left by the general FCC estimate    & 61               \\
            short directions blocked by the FCC core & 52               \\
            directions requiring the correction certificate      & 9                \\
            reduced-radius certificate: tested / accepted / depth & 9213 / 6912 / 11 \\
            minimum squared separation surplus & \(29/2500\) \\
            \bottomrule
        \end{tabular}
    \end{center}

    The transverse-basis convention below fixes the reported cell counts.
    A basis of integer vectors is \emph{unimodular} when the determinant of its
    basis matrix is \(\pm1\). In this setting, a unimodular
    transverse basis \((\vect{b}_1,\vect{b}_2)\) is one for which
    \((\vect{a},\vect{b}_1,\vect{b}_2)\) forms an integer basis of \(\Z^3\).
    A different unimodular transverse basis may give different cells and
    counts; the same test still gives a finite proof if it accepts every
    terminal cell.

    \subsection{Language-neutral exact pseudocode}

    Set \(\vect{C}_i=10\vect{c}_i\in\Z^3\) and
    \(U_i=100r_i\in\Z\), so the period is \(20\) in center coordinates.
    The pseudocode writes \texttt{C[i]} for \(\vect{C}_i\); vector symbols
    are unbolded in code. In the pseudocode, \texttt{dot}, \texttt{cross},
    \texttt{gcd}, \texttt{floor},
    and all divisions have their exact integer or rational meanings;
    \texttt{norm2(x)=dot(x,x)}.
    Loop bodies are indented and have explicit closing markers. A
    \texttt{repeat} loop runs until a \texttt{return} exits the procedure.
    In a one-line conditional, all statements after the colon belong to that
    conditional; \texttt{break} exits the innermost loop and \texttt{continue}
    starts its next iteration.

    \subsubsection{Exact packing check}
    The following computes the minimum squared separation surplus over all
    relevant pairs of periodic ball translates, as justified in
    Proposition~\ref{prop:packing-period}.
    \label{alg:packing}
    \begingroup
    \scriptsize
    \begin{verbatim}
PACKING-CHECK(C[1..38], U[1..38])
best := +infinity
for i := 1,...,38:
  for j := i,...,38:
    for k in {-1,0,1}^3:
      if i=j and k=(0,0,0): continue
      D := C[i] - C[j] - 20*k
      S := 100*norm2(D) - (U[i]+U[j])^2
      best := min(best,S)
    end for
  end for
end for
assert best = 116
output packing surplus best/10000 = 29/2500
    \end{verbatim}
    \endgroup

    \subsubsection{Periods and short \(\FCC\)-directions}
    \label{alg:period-directions}
    The bases and expected radii referenced below are the six rows of
    Table~\ref{tab:short-fcc-radii}.
    \begingroup
    \scriptsize
    \begin{verbatim}
PERIOD-AND-DIRECTION-CHECKS(C[1..38], U[1..38])
labelled := {(C[i] mod 20,U[i]) : i=1,...,38}
assert length(labelled)=38
for T in {(0,10,10),(10,0,10),(10,10,0)}:
  assert some (C[i]+T mod 20,U[i]) is not in labelled
end for
output full period lattice 2*Z^3

counts := empty map; correction := empty list
for a in {-3,-2,...,3}^3:
  if a=0 or gcd(|a[1]|,|a[2]|,|a[3]|) != 1: continue
  if first nonzero coordinate of a is negative: continue
  w := a if sum(a) is even, and 2*a otherwise
  if norm2(w) >= 16: continue
  kind := decreasing sort of (|w[1]|,|w[2]|,|w[3]|)
  counts[kind] := counts[kind]+1
  if kind in {(1,1,0),(2,0,0)}: append a to correction
end for
assert counts in the six table rows are (6,3,12,4,12,24), totaling 61
assert length(correction)=9

GRAM(u,v,w) := dot(u,v)-dot(u,w)*dot(v,w)/norm2(w)
for (w,b1,b2,expected) in the six table rows:
  assert |dot(w,cross(b1,b2))|=2
  (g11,g12,g22) := (GRAM(b1,b1,w),GRAM(b1,b2,w),GRAM(b2,b2,w))
  assert g12<=0 and g11+g12>=0 and g22+g12>=0
  mu2 := g11*g22*(g11+g22+2*g12)/(4*(g11*g22-g12^2))
  assert mu2=expected
end for
assert sum(counts[kind] for rows with expected<(7/10-1/250)^2)=52
return correction
    \end{verbatim}
    \endgroup

    \subsubsection{Adaptive line-space certificate}
    \label{alg:line-certificate}
    The input \texttt{correction} is the list of nine directions
    returned by the procedure in Appendix~\ref{alg:period-directions}.
    For each of the nine directions, fix a primitive
    integer representative \(\vect{a}\) and choose a transverse basis as follows.
    For an integer vector \(\vect{u}\), \texttt{linfnorm(u)} denotes
    \(\max_i\abs{u_i}\). Integer vectors are ordered first by increasing
    \texttt{linfnorm} and then lexicographically. Thus ``first'' below specifies
    a unique vector. The function \texttt{nearest} rounds an exact rational to a
    nearest integer, resolving a half-integer toward zero, and
    \[
        Q(\vect{u},\vect{v})
        =\norm{\vect{a}}^2\vect{u}\dotprod\vect{v}
        -(\vect{u}\dotprod\vect{a})(\vect{v}\dotprod\vect{a}).
    \]
    By~\eqref{eq:projected-inner-product},
    \[
        Q(\vect{u},\vect{v})
        =\norm{\vect{a}}^2
        \orthproj{a}(\vect{u})\dotprod\orthproj{a}(\vect{v}).
    \]
    Thus the loop in \texttt{TRANSVERSE-BASIS} applies
    the Gauss reduction of Lemma~\ref{lem:gauss-reduction} to
    \(\orthproj{a}(\vect{b}_1)\) and
    \(\orthproj{a}(\vect{b}_2)\).
    Subtracting multiples of \(\vect{a}\) from
    \(\vect{b}_1,\vect{b}_2\) leaves the projections unchanged.
    The certificate does not require a final sign change to make the
    projected inner product nonpositive.

    To choose the periodic copies, write each scaled center in the integer
    basis as
    \[
        \vect{C}_i=\alpha\vect{a}+\beta\vect{b}_1+\gamma\vect{b}_2,
        \qquad \alpha,\beta,\gamma\in\Z.
    \]
    Translations by \(20\vect{a}\) have zero projection onto
    \(\orthcomp{a}\). Set
    \[
        \vect{C}_0=\vect{C}_i
        -20\lfloor\beta/20\rfloor\vect{b}_1
        -20\lfloor\gamma/20\rfloor\vect{b}_2.
    \]
    This is the scaled center of a periodic copy of the same ball, with its
    coefficients of \(\vect{b}_1,\vect{b}_2\) in \([0,20)\).
    The procedure \texttt{NEARBY-COPIES} tests the \(25\) centers
    \[
        \vect{C}_0+20s\vect{b}_1+20t\vect{b}_2,
        \qquad s,t\in\{-2,-1,0,1,2\}.
    \]
    Each
    accepted cell has a specific translate satisfying the cell-acceptance test,
    so omitting other translates can only prevent acceptance, not create a
    false acceptance. Termination with every terminal cell accepted therefore
    proves that this window suffices.
    \begingroup
    \scriptsize
    \begin{verbatim}
TRANSVERSE-BASIS(a)
b1 := first nonzero integer vector with gcd(components(a cross b1))=1
normal := a cross b1
b2 := first integer vector with dot(normal,b2)=1
repeat:
  b1 := b1-nearest(dot(b1,a)/norm2(a))*a
  b2 := b2-nearest(dot(b2,a)/norm2(a))*a
  if Q(b2,b2)<Q(b1,b1): swap b1,b2
  m := nearest(Q(b1,b2)/Q(b1,b1))
  if m=0: return (b1,b2)
  b2 := b2-m*b1
end repeat

NEARBY-COPIES(a,b1,b2)
B := integer matrix with columns (a,b1,b2); copies := empty list
for each (C[i],U[i]):
  (alpha,beta,gamma) := B^(-1)*C[i]       // integers since |det B|=1
  C0 := C[i]-20*floor(beta/20)*b1-20*floor(gamma/20)*b2
  for s,t in {-2,-1,0,1,2}:
    append (C0+20*s*b1+20*t*b2,U[i]) to copies
  end for
end for
return copies

CERTIFY-DIRECTION(a)
(b1,b2) := TRANSVERSE-BASIS(a)
copies := NEARBY-COPIES(a,b1,b2)
L := l1norm(b1)+l1norm(b2); direction_norm2 := norm2(a)
stack := {(0,0,0)}; tested := accepted := deepest := 0
while stack is not empty:
  remove (i,j,d) from stack; tested := tested+1
  deepest := max(deepest,d); N := 2^d; covered := false
  X := 10*((2*i+1)*b1+(2*j+1)*b2)
  for (Z,U) in copies:
    D := X-N*Z
    distance_numerator := direction_norm2*norm2(D)-dot(D,a)^2
    margin_numerator := (5*U-2)*N-500*L
    if margin_numerator>0 and
       2500*distance_numerator < direction_norm2*margin_numerator^2:
      covered := true; break
    end if
  end for
  if covered: accepted := accepted+1; continue
  assert d<13
  add (2*i,2*j,d+1), (2*i+1,2*j,d+1),
      (2*i,2*j+1,d+1), (2*i+1,2*j+1,d+1) to stack
end while
return (tested,accepted,deepest)
    \end{verbatim}
    \noindent\begin{minipage}{\linewidth}
    \begin{verbatim}
CERTIFY-ALL-DIRECTIONS(correction)
total_tested := 0
total_accepted := 0
max_depth := 0
for a in correction:
  (tested, accepted, deepest) := CERTIFY-DIRECTION(a)
  total_tested := total_tested + tested
  total_accepted := total_accepted + accepted
  max_depth := max(max_depth, deepest)
end for
result := (total_tested, total_accepted, max_depth)
assert result = (9213, 6912, 11)
return result
    \end{verbatim}
    \end{minipage}
    \endgroup

\end{document}